\documentclass[11pt,a4paper]{amsart}
\usepackage[T1]{fontenc}
\usepackage{lmodern}
\usepackage[margin=1in]{geometry}
\usepackage{amsmath,amssymb,amsthm,mathtools,microtype,enumitem}
\usepackage{needspace}
\usepackage{xcolor}
\usepackage{tikz}
\usepackage[colorlinks=true,linkcolor=blue!45!black,citecolor=blue!45!black,urlcolor=blue!45!black]{hyperref}
\DeclareHookRule{begindocument}{microtype}{before}{amsart}
\setlist{itemsep=3pt,topsep=5pt,leftmargin=12pt}
\numberwithin{equation}{section}
\newtheorem{theorem}{Theorem}[section]
\newtheorem{proposition}[theorem]{Proposition}
\newtheorem{lemma}[theorem]{Lemma}
\newtheorem{corollary}[theorem]{Corollary}
\theoremstyle{remark}\newtheorem{remark}[theorem]{Remark}

\title[Strict hot spots on simply connected planar H\"older domains]{Strict hot spots and absence of interior critical points
on simply connected planar H\"older domains}
\author{Yutian Li}
\address{School of Mathematics and Statistics, Nanfang College, Guangzhou, 510970, Guangdong, China}
\email{yutianlee@gmail.com}
\date{30 November 2026}
\subjclass[2020]{Primary 35P05; Secondary 35B38, 35J25, 30C20.}
\keywords{Strict hot spots, Neumann eigenfunctions, interior critical points, H\"older domains, conformal maps.}
\hypersetup{
 pdftitle={Strict hot spots and absence of interior critical points on simply connected planar H\"older domains},
 pdfauthor={Yutian Li},
 pdfsubject={First positive Neumann eigenfunctions and global eigenfunction coordinates},
 pdfkeywords={Strict hot spots, Neumann eigenfunctions, H\"older domains, conformal maps}
}
\begin{document}
\begin{abstract}
Every nonzero real eigenfunction for the first positive Neumann
eigenvalue on a bounded simply connected planar H\"older domain has
no interior critical point. The H\"older condition is quasihyperbolic;
equivalently, the derivative of a Riemann map from the unit disk
is integrable to some power greater than two.
This includes Lipschitz domains and polygons with reentrant corners.
Multiple eigenvalues and non-Jordan boundaries are allowed.
The proof uses boundary primitives, a uniform multiplier estimate,
and approximation of the whole limiting eigenspace.
When this eigenspace is two-dimensional, every real basis gives global
coordinates onto a domain star-shaped about the origin.
A stream function proves injectivity and yields a Poisson formula
for the boundary angular derivative. Composing the coordinate map with
a Riemann map gives a boundary map with a Lipschitz inverse and a
rectifiable Jordan image. Its boundary angle is a circle diffeomorphism
with nonvanishing H\"older-continuous derivative. The image can be
nonconvex even when the source is smooth and strictly convex.
\end{abstract}
\maketitle

\tableofcontents
\clearpage

\section{Introduction and main results}
\label{sec:allp2-theorem}

The hot spots conjecture asks whether a first positive Neumann
eigenfunction attains its maximum and minimum on the boundary.
The question comes from the long-time behavior of heat flow with
insulating boundary conditions \cite[Section~1]{BanuelosBurdzy1999}.
Its strict form excludes all interior global extrema. We prove a
stronger conclusion on bounded simply connected planar H\"older
domains, defined by \eqref{eq:intro-qhb}: every such eigenfunction
has no interior critical point.
In particular, interior saddle points are also excluded.

The smooth-domain theorem of OpenAI \cite{OAI} gives this stronger
conclusion on smooth bounded simply connected planar domains, allowing
the first positive eigenvalue to be multiple. We extend the conclusion
to H\"older domains. This class includes all polygonal Jordan
domains, with reentrant corners allowed, and rough domains whose boundary
need not be Jordan. For the latter, strict hot spots is formulated using
the continuous pullback to the closed disk in
\eqref{eq:intro-prime-end-hotspots}.

Let $\Omega\subset\mathbb C$ be a bounded simply connected domain.
We define the Neumann Laplacian as the nonnegative self-adjoint
operator in $L^2(\Omega)$ associated with the quadratic form
$\int_\Omega|\nabla u|^2$ with form domain $H^1(\Omega)$.
Whenever the embedding $H^1(\Omega)\hookrightarrow L^2(\Omega)$
is compact, we list its eigenvalues with multiplicity as
$0=\mu_1<\mu_2\le\mu_3\le\cdots$.
In that case, for its first positive eigenvalue $\mu=\mu_2$, write
\[
 E_\mu=\left\{u\in H^1(\Omega;\mathbb R):
 \int_\Omega\nabla u\cdot\nabla\eta
   =\mu\int_\Omega u\eta
 \quad\hbox{for every }\eta\in H^1(\Omega)\right\}.
\]
The main results concern nonzero real eigenfunctions.

For $d_\Omega(x)=\operatorname{dist}(x,\partial\Omega)$ let
$k_\Omega$ denote quasihyperbolic distance: the infimum of
$\int_\gamma |d\xi|/d_\Omega(\xi)$ over rectifiable curves in
$\Omega$ joining the two points.
We call $\Omega$ a H\"older domain if, for some $x_0\in\Omega$ and
constants $C_1,C_2>0$,
\begin{equation}\label{eq:intro-qhb}
 k_\Omega(x_0,x)\le
 C_1\log\frac{d_\Omega(x_0)}{d_\Omega(x)}+C_2
 \qquad(x\in\Omega).
\end{equation}
Thus the H\"older condition used here concerns quasihyperbolic
distance and does not require the boundary to be a local graph.
We write $\mathbb D=\{z\in\mathbb C:|z|<1\}$,
$\mathbb S^1=\partial\mathbb D=\{z\in\mathbb C:|z|=1\}$, and $\delta(z)=1-|z|$.
The measure $dA$ is planar Lebesgue measure. For $0<s<\infty$,
$A^s$ denotes the Bergman space of holomorphic functions on
$\mathbb D$ belonging to $L^s(\mathbb D,dA)$.

\begin{theorem}[Absence of interior critical points and strict hot spots]\label{thm:allp2-main}
Let $\Omega$ be a bounded simply connected planar H\"older domain.
Then the Neumann Laplacian with form domain $H^1(\Omega)$ has
compact resolvent. Denote its first positive eigenvalue by
$\mu=\mu_2(\Omega)$. Every nonzero $u\in E_\mu$ has no interior
critical point:
\begin{equation}\label{eq:intro-main}
 \nabla u(x)\ne0
 \qquad(x\in\Omega).
\end{equation}
Consequently, the strict hot spots property holds: if
$\Phi:\mathbb D\to\Omega$ is a conformal bijection, then
$v=u\circ\Phi$ extends continuously to $\overline{\mathbb D}$ and
\begin{equation}\label{eq:intro-prime-end-hotspots}
 \min_{\mathbb S^1}v<v(z)<\max_{\mathbb S^1}v
 \qquad(z\in\mathbb D).
\end{equation}
Neither simplicity of the eigenvalue nor a Jordan boundary is required.
\end{theorem}

For a conformal bijection $\Phi:\mathbb D\to\Omega$,
Proposition~\ref{prop:qhb-conformal-regularity} identifies the
H\"older-domain hypothesis with $\Phi'\in A^{2p}$ for some $p>1$.
This integrability implies boundedness of $\Omega$ and compactness
of $H^1(\Omega)\hookrightarrow L^2(\Omega)$ by
Lemma~\ref{lem:sa-conformal}; the latter gives the compact-resolvent
conclusion of Theorem~\ref{thm:allp2-main}.

For a Jordan domain, these conformal boundary values are the usual
boundary values of $u$. For a non-Jordan domain,
\eqref{eq:intro-prime-end-hotspots} is stated on the conformal boundary;
distinct approaches to the same Euclidean boundary point may carry
different values. In particular, the two banks of a slit correspond to
distinct boundary approaches.

Theorem~\ref{thm:allp2-main} applies to every bounded simply connected
planar Lipschitz domain, including every polygonal Jordan domain with
reentrant corners. It therefore answers Question~1 in
\cite[Section~6]{Banuelos2026}, which asks for non-strict hot spots on
this class without assuming simplicity, and proves the stronger
absence of all interior critical points.
It also applies to all bounded convex planar
domains, bounded simply connected John domains, and quasidisks.
Section~\ref{sec:qhb-scope} defines John domains and quasidisks and
explains these inclusions. Section~\ref{subsec:domain-examples} gives
the von Koch snowflake and a non-Jordan slit domain as concrete examples.

When $E_\mu$ is two-dimensional, its eigenfunctions give global
coordinates with the following boundary behavior.

\begin{corollary}[Global eigenfunction coordinates]\label{cor:intro-eigen-coordinates}
Let $\Omega$ satisfy the hypotheses of Theorem~\ref{thm:allp2-main}
and suppose $\dim_{\mathbb R}E_\mu=2$. For any real basis $u_1,u_2$
of $E_\mu$, define
\[
 \mathcal U(x)=(u_1(x),u_2(x)),\qquad \Sigma=\mathcal U(\Omega).
\]
Then $\mathcal U$ is a real-analytic diffeomorphism onto the bounded Jordan
domain $\Sigma$, which is star-shaped about the origin. The
eigenfunctions $u_1,u_2$ have exactly one common zero in $\Omega$.
For every conformal bijection $\Phi:\mathbb D\to\Omega$, the map
$\mathcal U\circ\Phi$ extends to a homeomorphism from $\overline{\mathbb D}$
onto $\overline\Sigma$, and the inverse boundary map
$\partial\Sigma\to\mathbb S^1$ is Lipschitz. Every ray from the origin meets
$\partial\Sigma$ exactly once, at a distance that is a continuous
positive function of its angle. The boundary $\partial\Sigma$ is
rectifiable, and for some $\beta>0$ the map
$\zeta\mapsto(\mathcal U\circ\Phi)(\zeta)/|(\mathcal U\circ\Phi)(\zeta)|$
is a $C^{1,\beta}$ diffeomorphism from $\mathbb S^1$ onto the unit
circle in $\mathbb R^2$.

For every nonzero $w\in E_\mu$, the boundary trace of $w\circ\Phi$
has exactly two zeros on $\mathbb S^1$, both simple in the conformal
parameter. The closed nodal set of $w\circ\Phi$ in
$\overline{\mathbb D}$ is a single arc, $C^{1,\beta}$ up to its
two endpoints, and meets $\mathbb S^1$ orthogonally.
If $\Omega$ is Jordan, $\mathcal U$ extends
to a homeomorphism from $\overline\Omega$ onto $\overline\Sigma$.
\end{corollary}

The proof is given in Corollary~\ref{cor:local-eigen-coordinates}
and Proposition~\ref{prop:eigen-coordinate-boundary}, with the nodal
and radial regularity in Corollaries~\ref{cor:eigen-boundary-zeros}
and~\ref{cor:eigen-radial-regularity}. These statements give explicit
regularity exponents in terms of the integrability of $\Phi'$.

\subsection{Background and earlier results}
\label{subsec:hot-spots-background}

Rauch raised the question at Tulane in 1974; see the account in
\cite[Section~1]{BanuelosBurdzy1999}. On a smooth bounded connected
domain, the leading decaying term of the Neumann heat flow is the
projection of the initial temperature onto the first positive
eigenspace, multiplied by its exponential decay factor, provided this
projection is nonzero \cite[Proposition~2.1]{BanuelosBurdzy1999}.
Thus the initial temperature selects an eigenfunction in that
eigenspace, which explains the need to treat multiplicity.

For an eigenfunction continuous on the closure, the non-strict hot
spots property places every interior value between the minimum and
maximum boundary values. The strict property requires both inequalities
to be strict at every interior point. One must also distinguish a
statement for every first positive eigenfunction from the existence
of one eigenfunction with the property. These are separate formulations
in \cite[Section~1]{BanuelosBurdzy1999}.

The unrestricted planar conjecture is false. Burdzy and Werner
\cite{BurdzyWerner1999} constructed a bounded planar domain with two
holes and a simple first positive Neumann eigenvalue whose eigenfunction
has a strict interior maximum. Burdzy \cite[Theorem~1.1]{Burdzy2005}
gave a domain with one hole for which both extrema lie strictly beyond
the range of boundary values. The simply connected planar version is
stated in \cite[Conjecture~1.2(ii)]{Burdzy2005}.
There is also an obstruction in higher dimensions: de Dios Pont
\cite{deDiosPont2024} constructed smooth convex counterexamples in
sufficiently large dimensions.

Jerison and Nadirashvili
\cite[Theorem~1.4]{JerisonNadirashvili2000} proved strict hot spots on
bounded convex planar domains with two axes of symmetry, including
multiple eigenvalues. Atar and Burdzy \cite[Theorem~1]{AtarBurdzy2004}
proved strict hot spots for lip domains: bounded Lipschitz planar
domains enclosed by two graphs with Lipschitz constants at most one,
meeting at both endpoints of their common interval.
This is a special class of Lipschitz domains.
Judge and Mondal proved absence of interior
critical points on triangles; their original paper should be read
together with its correction \cite{JudgeMondal2020,JudgeMondal2022}.
Chen, Gui, and Yao obtained finer conclusions on nonvertex critical
points and directional monotonicity \cite{ChenGuiYao2026}, while
Hatcher treated L-shaped domains and other nonconvex polygonal families
\cite{Hatcher2025}.
Deng, Gui, Jiang, Yang, Yao, and Zou
\cite[Theorem~1.3]{DengGuiJiangYangYaoZou2026} prove absence of
interior critical points on isosceles trapezoids and kites, including
the cases of a two-dimensional first positive eigenspace. They also
classify symmetry transitions and boundary critical points on these
domains.

Burdzy \cite{Burdzy2026} derives a non-strict hot-spots conclusion
on bounded simply connected domains with compact Neumann resolvent
and simple first positive eigenvalue. Ba\~nuelos \cite{Banuelos2026}
treats finite-area simply connected domains under compactness
and simplicity, with extrema approached in arbitrarily thin
boundary layers. These results use the smooth theorem and apply
to geometric classes that are not all covered here.
Deng, Jiang, and Yang \cite{DengJiangYang2026} prove strict boundary
maximum modulus on every bounded convex planar domain without
smoothness or simplicity assumptions. Theorem~\ref{thm:allp2-main}
excludes all interior critical points on that class; their results
for higher eigenvalues concern different questions.

Pineda Montoya, Gallego Gonz\'alez, Porras Barrera, and R\'ua Mu\~noz
\cite{Pineda2026} prove gradient nonvanishing on piecewise smooth
Jordan domains whose corner angles lie in $(0,\pi)$, allowing
multiplicity. Their domains need not be globally convex. The
present theorem permits reentrant corners and non-Jordan boundary
under the quasihyperbolic condition.

Rohleder's variational approach identifies gradients of Neumann
eigenfunctions through a vector-valued minimization principle
\cite{Rohleder2023}. From the smooth-domain proof in \cite{OAI},
we use the reciprocal Green-kernel method, its scalar multiplier
identity, and its argument using the support of the boundary trace.
The realization by squared Hilbert-space distances in
Proposition~\ref{prop:rkchecked-realization} and the Gaussian kernels
in \eqref{eq:coherent-Gram} use Schoenberg's construction
\cite{Schoenberg1938}. Appendix~\ref{sec:rkchecked} proves the reciprocal-kernel theorem
using the regularization of \cite[Section 5]{OAI} and the global
matrix flow of \cite[Lemmas 4.4--4.5]{Pineda2026}, with the limiting
steps required here.
We also construct a boundary primitive for the Neumann gradient and
control multiplied fields with small positive energy. Approximation
of the whole first positive eigenspace then handles eigenvalue splitting.

Eigenfunctions have also been used to construct local coordinates.
Jones, Maggioni, and Schul \cite[Theorem~2.1.1]{JonesMaggioniSchul2010} obtain local
parametrizations with quantitative distortion bounds by choosing
appropriate Laplace eigenfunctions. Corollary~\ref{cor:intro-eigen-coordinates}
applies globally to every real basis of the first positive eigenspace
when its dimension is two. The nonvanishing Jacobian gives local
invertibility. A positive stream function then verifies the global
inversion criterion of Gorni and Zampieri
\cite[Theorem~1.1 and Proposition~1.2]{GorniZampieri1994}, which proves
injectivity and star-shapedness. Boundary elliptic regularity and a
Poisson formula make the angular derivative strictly positive after
an orientation is chosen. This gives the boundary homeomorphism,
the Lipschitz inverse, and the regularity of nodal endpoints.

Star-shapedness does not imply convexity of the image. The
counterexamples of B\'erard and Helffer
\cite[Theorem~1.2 and Proposition~3.5]{BerardHelffer2021} give smooth
strictly convex source domains with nonconvex eigenfunction-coordinate
images. The equilateral triangle gives the deltoid map described by
Bakry and Zribi \cite[Section~3]{BakryZribi2016}. We use this established
map to exhibit the loss of radial and angular regularity at a corner
and an image that is star-shaped only about the origin; see
Section~\ref{subsec:eigen-coordinate-convexity}.

\subsection{The quantitative multiplier estimate}

Two difficulties arise when a smooth exhaustion approaches a rough
domain. Nonvanishing of the gradient is not preserved by
$C^1_{\mathrm{loc}}$ convergence, and a multiple first positive
eigenvalue may split into several nearby eigenvalues. We address both
difficulties by estimating the energy of multiplied vector fields and
approximating the whole limiting eigenspace.

We state the estimate on the smooth approximating domains.
Let $\Phi:\mathbb D\to\Omega$ be a conformal bijection, suppose
$\Phi'\in A^{2p}$ for some $p>1$, and put $r=2p/(p+1)>1$.
For $\rho<1$ close to one, set
\[
 \Phi_\rho(z)=\Phi(\rho z),\quad
 \Omega_\rho=\Phi(\rho\mathbb D),\quad
 \mu_\rho=\mu_2(\Omega_\rho),\quad
 q_\rho=\mu_\rho|\Phi_\rho'|^2.
\]
Let $A$ denote ordinary Poisson extension on the disk, with respect to
circle arclength, and let
$G=(-\Delta_{\mathbb D,D})^{-1}$ be the Dirichlet Green operator.
On $L^2(\mathbb D,q_\rho\,dA)$, let $Gq_\rho$ denote
$f\mapsto G(q_\rho f)$, and let $I$ be the identity operator.
For $g\in H^{1/2}(\mathbb S^1;\mathbb C)$, set
\[
 S_\rho=(I-Gq_\rho)^{-1},\qquad V_{\rho,g}=S_\rho Ag.
\]
The inverse is justified by Lemma~\ref{lem:sa-resolvent}.
We identify a planar vector $(X_1,X_2)$ with $X_1+iX_2$ and use
the Wirtinger derivatives
$\partial_z=(\partial_x-i\partial_y)/2$ and
$\partial_{\bar z}=(\partial_x+i\partial_y)/2$.
For $V\in H^1(\mathbb D;\mathbb C)$, define
\[
 \mathcal E_\rho(V)=
 \int_{\mathbb D}\bigl(4|V_z|^2-q_\rho|V|^2\bigr)\,dA.
\]
For a real Hilbert-valued multiplier, extension is taken coordinatewise
and the energies of the resulting complex fields are summed. Lemma~\ref{lem:sa-hodge} and
\eqref{eq:sa-conformal-form} imply that it is nonnegative when the
corresponding vector field is tangent to $\partial\Omega_\rho$.

Fix $z_0\in\mathbb D$. Let $K_\rho=K_{q_\rho}$ be the Poisson kernel
for $-\Delta-q_\rho$ on $\mathbb D$, with respect to circle arclength,
and let $N_\rho(s,t)=-\partial_{\nu_s}K_\rho(s,t)$ for distinct
$s,t\in\mathbb S^1$,
where $\nu_s=s$ is the outward unit normal and the derivative is taken
in the first variable.
The chain formulas \eqref{eq:hk-K}--\eqref{eq:hk-R} give these kernels
and their normalization. Proposition~\ref{prop:rkchecked-realization}
constructs a map $b_\rho$ from $\mathbb S^1$ into a real Hilbert space
$\mathcal H_\rho$ such that
\[
 \|b_\rho(s)-b_\rho(t)\|^2
 =\frac{K_\rho(z_0,s)K_\rho(z_0,t)}{N_\rho(s,t)}
 \qquad(s,t\in\mathbb S^1,\ s\ne t).
\]
We normalize $b_\rho$ to vanish at a fixed point of $\mathbb S^1$.
Proposition~\ref{prop:hk-curve} proves its $C^{1,\beta}$ regularity
for $0<\beta<\min\{1/2,1-1/p\}$.

\begin{theorem}[Uniform multiplier estimate]\label{thm:intro-multiplier}
With the preceding notation, there exist $\rho_0\in(0,1)$ and
$C>0$, depending only on $\Phi$, $p$, and $z_0$, with the following
property. Let $\rho\in(\rho_0,1)$ and let
$g\in C^\infty(\mathbb S^1;\mathbb C)$ satisfy
$g(s)/(is\Phi_\rho'(s))\in\mathbb R$ for every $s\in\mathbb S^1$.
This is the tangency condition on $\partial\Omega_\rho$. With
$\theta$ denoting the angular coordinate on $\mathbb S^1$, suppose,
in the sense of distributions, that
\[
 g=g_0+\partial_\theta g_1,\qquad
 g_0,g_1\in L^r(\mathbb S^1;\mathbb C),
 \qquad M=\|g_0\|_{L^r(\mathbb S^1)}+\|g_1\|_{L^r(\mathbb S^1)}.
\]
For $0<t\le1$, with the exponential interpreted pointwise,
\begin{equation}\label{eq:intro-multiplier}
 \mathcal E_\rho(V_{\rho,\exp(-t\|b_\rho\|^2/2)b_\rho g})
 \le C\left(
 \frac{\mathcal E_\rho(V_{\rho,g})}{t}
 +|V_{\rho,g}(z_0)|^2+tM^2\right).
\end{equation}
The constant $C$ is independent of $\rho$, $g$, $t$, and the
Hilbert dimension.
\end{theorem}

The proof is given in Theorem~\ref{thm:allp2-gaussian}.
Although the estimate uses a distributional representation of $g$,
the trace in the statement is smooth. An arbitrary distribution of
the displayed form need not have a finite-energy extension. The
limiting fields on the rough domain are obtained by spectral
compactness in Section~\ref{sec:allp2-spectral}.

For a Neumann gradient, the eigenfunction equation supplies such
representations with a uniform $L^r$ bound. If the gradient vanishes
at $\Phi(z_0)$, its approximations have energies $e_\rho\to0$ and
values at $z_0$ tending to zero. Choosing $t_\rho\downarrow0$ with
$e_\rho/t_\rho\to0$ makes the right side of
\eqref{eq:intro-multiplier} tend to zero. Approximation of the whole
limiting eigenspace then shows that the multiplied limits are again
eigenfunction gradients.

\subsection{Outline of the proof}

Let $\lambda_1^D(\Omega)$ denote the first Dirichlet eigenvalue and
put $q=\mu|\Phi'|^2$. The strict inequality
$\mu<\lambda_1^D(\Omega)$ makes $(I-Gq)^{-1}$ positivity preserving
and bounded.
Sections~\ref{sec:sa-spectral} and~\ref{sec:sa-hodge} prove the
spectral and weak divergence--curl statements, including
$\dim_{\mathbb R}E_\mu\le2$, without using reciprocal kernels.
Section~\ref{sec:holder-kernel-audit} proves boundary-kernel
regularity and convergence under strong $L^p$ convergence of the
potential, using the Green-function estimates of
ter Elst--Ouhabaz \cite{EO2019}.

Section~\ref{sec:primitive-trace-all-p2} constructs the boundary
primitive using the reciprocal-derivative estimate in
\eqref{eq:allp2-inverse}. This follows from
the non-strict integral-means bound proved in
\cite[Section~2.3, equation~(27)]{BIPS2026}; no strict improvement
of that bound is used. Sections~\ref{sec:allp2-spectral} and~\ref{sec:coherent-near-null}
combine spectral approximation with the multiplier estimate.
Section~\ref{sec:allp2-support} then derives a contradiction from the
support of the harmonic boundary trace, using Sawyer's
unique-continuation theorem \cite{Sawyer1984}.
Section~\ref{subsec:global-eigen-coordinates} constructs the stream
function and applies the global inversion criterion. Section~\ref{subsec:eigen-coordinate-boundary}
proves the boundary angular formula, uses boundary Schauder estimates
to obtain its regularity, and derives the boundary and nodal
consequences. Section~\ref{subsec:eigen-coordinate-convexity} treats
convexity, points about which the image is star-shaped, and the
equilateral triangle.
The characterization of H\"older domains in
Proposition~\ref{prop:qhb-conformal-regularity} uses
Smith--Stegenga \cite{SmithStegenga1991}.

The endpoint $\Phi'\in A^2$ is not covered. Constants in the
boundary-primitive estimates may deteriorate as $p\downarrow1$.
The argument does not cover all bounded simply connected domains
with compact Neumann resolvent, or the broader finite-area class in
\cite[Section~6, Question~2]{Banuelos2026}.

\section{H\"older domains and conformal integrability}
\label{sec:qhb-scope}

We first relate the quasihyperbolic boundary condition
\eqref{eq:intro-qhb} to integrability of the conformal derivative.

A bounded domain $\Omega$ is a John domain if there are $x_0\in\Omega$
and $c>0$ such that every $x\in\Omega$ can be joined to $x_0$ by an
arclength parametrized curve $\gamma:[0,L]\to\Omega$ satisfying
$\gamma(0)=x$, $\gamma(L)=x_0$, and
$d_\Omega(\gamma(s))\ge cs$ for $0\le s\le L$.

\begin{proposition}\label{prop:qhb-conformal-regularity}
Let $\Omega$ be a bounded simply connected planar domain and let
$\Phi:\mathbb D\to\Omega$ be a conformal bijection. The following
are equivalent:
\begin{enumerate}
\item $\Phi'\in L^{2p}(\mathbb D)$ for some $p>1$;
\item $\Omega$ satisfies the quasihyperbolic boundary condition
      \eqref{eq:intro-qhb}.
\end{enumerate}
Every bounded simply connected planar John domain satisfies these conditions.
\end{proposition}

\begin{proof}
Precomposing with a fixed disk automorphism preserves membership of
$\Phi'$ in $L^{2p}(\mathbb D)$, so normalize $\Phi(0)=x_0$.
Let $\lambda_\Omega$ and $h_\Omega$ denote hyperbolic density and
distance, normalized to curvature $-1$. Monotonicity of the density
and the Koebe quarter theorem give
\begin{equation}\label{eq:qhb-hyperbolic-comparison}
 \frac{1}{2d_\Omega(x)}\le\lambda_\Omega(x)
 \le\frac{2}{d_\Omega(x)},
 \qquad
 \frac12 k_\Omega(x,y)\le h_\Omega(x,y)\le2k_\Omega(x,y).
\end{equation}

Suppose first that \eqref{eq:intro-qhb} holds.  The theorem of
Smith--Stegenga \cite[Theorem A and Theorem 2]{SmithStegenga1991} gives
some $c_0>0$ such that
\begin{equation}\label{eq:qhb-exp-integrability}
 \int_\Omega e^{c_0 k_\Omega(x_0,x)}\,dA(x)<\infty.
\end{equation}
The theorem requires no Jordan, rectifiability, or boundary smoothness
assumption. Fix $0<\varepsilon\le c_0/2$. For $x=\Phi(z)$, conformal
invariance of hyperbolic distance and
\eqref{eq:qhb-hyperbolic-comparison} imply
\[
 \delta(z)^{-\varepsilon}
 \le\left(\frac{1+|z|}{1-|z|}\right)^\varepsilon
 =e^{\varepsilon h_\Omega(x_0,x)}
 \le e^{2\varepsilon k_\Omega(x_0,x)}.
\]
Since $\Phi$ is bounded, Cauchy's derivative estimate gives
$|\Phi'(z)|\le C\delta(z)^{-1}$. The conformal change of variables
therefore yields
\begin{align*}
 \int_{\mathbb D}|\Phi'|^{2+\varepsilon}\,dA
 &\le C^\varepsilon\int_{\mathbb D}\delta^{-\varepsilon}|\Phi'|^2\,dA\\
 &=C^\varepsilon\int_\Omega
       \delta(\Phi^{-1}(x))^{-\varepsilon}\,dA(x)\\
 &\le C^\varepsilon\int_\Omega e^{2\varepsilon k_\Omega(x_0,x)}\,dA(x)
 <\infty.
\end{align*}

Conversely, if $\Phi'\in L^{2p}(\mathbb D)$ for some $p>1$, then
Morrey's inequality on the disk gives a continuous extension
$\Phi\in C^{0,\beta}(\overline{\mathbb D})$, where
$\beta=1-1/p>0$.  Its values on the unit circle lie in
$\partial\Omega$. For $z\ne0$, comparison with the radial endpoint gives
\[
 d_\Omega(\Phi(z))\le |\Phi(z)-\Phi(z/|z|)|
 \le C\delta(z)^\beta.
\]
It follows that
\[
 h_\Omega(x_0,\Phi(z))
 =\log\frac{1+|z|}{1-|z|}
 \le\log2+\frac1\beta
       \log\frac{C}{d_\Omega(\Phi(z))}.
\]
Using $k_\Omega\le2h_\Omega$ proves
\eqref{eq:intro-qhb} after adjusting the additive constant.

For a John domain, take a curve as in the definition above. At the endpoint,
$L\le d_\Omega(x_0)/c$.  For $s\le d_\Omega(x)/2$, the
one-Lipschitz property of $d_\Omega$ gives
$d_\Omega(\gamma(s))\ge d_\Omega(x)/2$.  Splitting the curve integral
there gives
\[
 k_\Omega(x,x_0)
 \le 1+\frac1c
 \max\left\{0,\log\frac{2d_\Omega(x_0)}{c\,d_\Omega(x)}\right\},
\]
which implies \eqref{eq:intro-qhb}.
\end{proof}

\subsection{Examples}\label{subsec:domain-examples}

All domains in the following examples are bounded, planar, and
simply connected.

\begin{itemize}
\item{Lipschitz domains and polygons.}
A domain is Lipschitz if, after a rigid change of coordinates near
each boundary point, its boundary is the graph of a Lipschitz
function and the domain lies on one side of that graph.
Every bounded Lipschitz domain is a John domain
\cite[p.~582, following Definition~3.1]{BuckleyKoskela1995}.
Proposition~\ref{prop:qhb-conformal-regularity} therefore applies.
This includes all bounded simply connected $C^1$ domains and all
polygonal Jordan domains. Reentrant corners are allowed, so no
convexity assumption is needed for polygons.

\item{Convex domains.}
Every bounded convex domain is a John domain. Indeed, choose
$B(x_0,r_0)\subset\Omega$. Convexity gives
$B((1-t)x+tx_0,tr_0)\subset\Omega$ for $x\in\Omega$ and
$0<t\le1$. The segment from $x$ to $x_0$, parametrized by
arclength $s$, thus stays at distance at least
$r_0s/\operatorname{diam}\Omega$ from the boundary.

\item{Quasidisks and a fractal example.}
A quasidisk is the image of the disk under a quasiconformal
homeomorphism of the plane. A Riemann map $\Phi:\mathbb D\to\Omega$
onto a bounded quasidisk extends quasiconformally across the circle by reflection
\cite{Ahlfors1963}. Local higher integrability of the derivative of
this extension \cite[Corollary~2.2]{Astala1998} gives
$\Phi'\in L^{2+\eta}$ for some $\eta>0$. The exponent may depend
on the domain; no fixed bound on the quasiconformal distortion is
assumed. The domain enclosed by the von Koch snowflake gives a
concrete example with fractal boundary: it is a John domain
\cite[p.~578]{BuckleyKoskela1995}.

\item{A non-Jordan domain.}
Let $\Omega=\mathbb D\setminus[0,1)$ and
$D_+=\{w:|w|<1,\ \operatorname{Im}w>0\}$.
The domain $D_+$ is bounded and convex, so a conformal map
$\Psi:\mathbb D\to D_+$ satisfies $\Psi'\in A^{2p}$ for some
$p>1$ by Proposition~\ref{prop:qhb-conformal-regularity}.
Squaring maps $D_+$ conformally onto $\Omega$. Thus
$\Phi=\Psi^2$ is a conformal bijection onto $\Omega$ and
$|\Phi'|\le2|\Psi'|$. Theorem~\ref{thm:allp2-main} applies,
with the two banks of the slit represented by distinct conformal
boundary approaches.
\end{itemize}

\section{The first positive Neumann eigenspace}
\label{sec:sa-spectral}

Throughout this section, let $\Phi:\mathbb D\to\Omega$ be a
conformal bijection and suppose $F=\Phi'\in L^{2p}(\mathbb D)$
for some $p>1$. Put $J=|F|^2$.
The Dirichlet Laplacian is defined by the same quadratic form as the
Neumann Laplacian, with domain $H^1_0(\Omega)$ in place of $H^1(\Omega)$.
Write $0<\lambda_1^D(\Omega)\le\lambda_2^D(\Omega)\le\cdots$ for its
eigenvalues, repeated with multiplicity; compactness is proved in
Lemma~\ref{lem:sa-conformal}.
The arguments below use neither a boundary trace on $\Omega$ nor a
homeomorphic extension of $\Phi$ to the boundary.

\subsection{Sobolev spaces and regularity of eigenfunctions}

\begin{lemma}\label{lem:sa-conformal}
Under the conformal hypothesis of this section, $\Omega$ is bounded.
Composition with $\Phi$ is a bounded
isomorphism between $H^1(\Omega)$ and $H^1(\mathbb D)$, and between
$H^1_0(\Omega)$ and $H^1_0(\mathbb D)$. The variational Neumann and Dirichlet
resolvents are compact. If $u$ is a real Neumann eigenfunction at
eigenvalue $\lambda$, then $v=u\circ\Phi$ solves
\begin{equation}\label{eq:sa-disk-equation}
 \int_{\mathbb D}\nabla v\cdot\nabla\eta
   =\lambda\int_{\mathbb D}Jv\eta
 \qquad(\eta\in H^1(\mathbb D)),
\end{equation}
and $v\in W^{2,p}(\mathbb D)\cap C(\overline{\mathbb D})$.
\end{lemma}
\begin{proof}
The submean inequality for holomorphic functions gives
\[
 |F(z)|\le C_p\|F\|_{2p}\delta(z)^{-1/p}.
\]
Since $1/p<1$, radial integration bounds $\Phi$ uniformly.

The interior change of variables preserves Dirichlet energy:
\begin{equation}\label{eq:sa-conformal-energy}
 \int_\Omega|\nabla u|^2
   =\int_{\mathbb D}|\nabla(u\circ\Phi)|^2.
\end{equation}
For $u\in H^1(\Omega)$ the pullback is locally $H^1$ and has finite
global gradient norm. Fix an interior disk $B\Subset\mathbb D$.
The function $|F|$ is bounded above and away from zero on $B$, so
the $L^2(\Omega)$ norm controls the pullback's $L^2(B)$ norm.
The Poincar\'e inequality with an interior $L^2$ term,
\[
 \|w\|_{L^2(\mathbb D)}
 \le C_B\bigl(\|\nabla w\|_{L^2(\mathbb D)}+\|w\|_{L^2(B)}\bigr)
\]
applied to bounded truncations, followed by Fatou's lemma, puts the
pullback in $H^1(\mathbb D)$ and gives a bounded composition operator.
For the converse, let $p'=p/(p-1)$. H\"older's inequality and
Sobolev embedding on the disk give
\[
 \int_\Omega |v\circ\Phi^{-1}|^2
 =\int_{\mathbb D}J|v|^2
 \le\|J\|_p\|v\|_{2p'}^2
 \le C\|v\|_{H^1(\mathbb D)}^2.
\]
Together with \eqref{eq:sa-conformal-energy}, this proves the
$H^1$ isomorphism. Smooth compactly supported functions transform
to smooth compactly supported functions; density and the same bounds
give the $H^1_0$ correspondence.

The compact disk embedding $H^1(\mathbb D)\Subset L^{2p'}(\mathbb D)$
gives compactness in $L^2(\mathbb D,J\,dA)$, and hence compactness of
$H^1(\Omega)\hookrightarrow L^2(\Omega)$. Since $\Omega$ is connected,
the Neumann kernel consists of constants and $\mu_2(\Omega)>0$.
The Dirichlet assertion follows either by restriction or by zero
extension to a containing rectangle.

The correspondence of the two $H^1$ spaces proves
\eqref{eq:sa-disk-equation}. Sobolev embedding gives $v\in L^s$ for
every finite $s$. Choose $1<p_0<p$; H\"older's inequality gives
$Jv\in L^{p_0}$.
The Neumann estimate on the smooth disk gives $v\in W^{2,p_0}$,
hence $v\in L^\infty$. Now $Jv\in L^p$, and the same estimate
gives $v\in W^{2,p}\subset C(\overline{\mathbb D})$.
\end{proof}

\subsection{The strict Dirichlet gap}

\begin{lemma}\label{lem:sa-gap}
With $\mu=\mu_2(\Omega)$, one has
\begin{equation}\label{eq:sa-gap}
 0<\mu<\lambda_1^D(\Omega).
\end{equation}
\end{lemma}
\begin{proof}
Put $\lambda=\lambda_1^D(\Omega)$ and choose a positive Dirichlet
ground state $\phi_1\in H^1_0(\Omega)$. For every $k\in\mathbb R^2$
with $|k|^2=\lambda$, let
\[
 w_k=e^{ik\cdot x}
       -\frac{\int_\Omega e^{ik\cdot x}}{\int_\Omega\phi_1}\,\phi_1.
\]
It has mean zero and is nonzero. Define
\[
 Q(f,h)=\int_\Omega\nabla f\cdot\overline{\nabla h}
             -\lambda\int_\Omega f\bar h.
\]
The eigenfunction equation gives $Q(\phi_1,\phi_1)=0$, and
$Q(e^{ik\cdot x},e^{ik\cdot x})=0$ pointwise. Testing the
plane-wave equation against $\phi_1\in H^1_0$ gives
$Q(e^{ik\cdot x},\phi_1)=0$. Thus $Q(w_k,w_k)=0$, and the
Neumann variational principle gives $\mu\le\lambda$.

If equality held, every $w_k$ would be a minimizer on the
mean-zero subspace and hence a member of the complexification of
$E_\mu$. All plane waves with $|k|^2=\lambda$ would lie in the
complex linear span of $E_\mu$ and $\phi_1$, a finite-dimensional
space. But arbitrarily large finite collections of distinct plane
waves are linearly independent on every nonempty open set: restrict
an analytic relation to a line whose direction separates their
frequencies and use the Vandermonde determinant. This is a contradiction.
\end{proof}

\subsection{A positive resolvent and a boundary functional}

Set $q=\mu J$ and let $G=(-\Delta_{\mathbb D,D})^{-1}$.
Use the real or complex weighted space
$L^2(q)=L^2(\mathbb D,q\,dA)$. Define $T=G(q\,\cdot)$.
The gap and the $H^1_0$ correspondence give
\begin{equation}\label{eq:sa-form-bound}
 \int_{\mathbb D}q|w|^2\le d\int_{\mathbb D}|\nabla w|^2,
 \qquad d=\frac{\mu}{\lambda_1^D(\Omega)}<1,
 \qquad w\in H^1_0(\mathbb D).
\end{equation}

\begin{lemma}\label{lem:sa-resolvent}
The operator $T$ on $L^2(q)$ is self-adjoint, nonnegative, positivity
preserving, and has norm at most $d$. Hence
$S=(I-T)^{-1}=\sum_{j\ge0}T^j$ is bounded and positivity preserving.
Let $A$ denote ordinary disk Poisson extension. Then
\begin{equation}\label{eq:sa-positive-functional}
 L(h):=\int_{\mathbb D}q\,SAh\,dA,\qquad h\in C(\mathbb S^1;\mathbb R),
\end{equation}
defines a continuous real linear functional satisfying
\[
 h\ge0,\ h\not\equiv0\quad\Longrightarrow\quad L(h)>0.
\]
For every $u\in E_\mu$, with $v=u\circ\Phi$ and
$h=v|_{\mathbb S^1}$, one has $v=SAh$ and $L(h)=0$.
\end{lemma}
\begin{proof}
Let $B:H^1_0(\mathbb D)\to L^2(q)$ be the inclusion, with the
Dirichlet inner product on $H^1_0$. The form bound gives
$\|B\|\le\sqrt d$, and $T=BB^*$ proves self-adjointness,
nonnegativity, and $\|T\|\le d$. Positivity preservation follows
from the positive Dirichlet Green function. The series for $S$
therefore converges in operator norm and preserves positivity.

Since $\|Ah\|_{L^2(q)}\le\|q\|_1^{1/2}\|h\|_\infty$,
\[
 |L(h)|\le\frac{\|q\|_1}{1-d}\|h\|_\infty.
\]
If $h\ge0$ is continuous and nonzero, the ordinary Poisson extension
is strictly positive throughout the disk. Since $SAh\ge Ah$ and
$q=\mu|F|^2>0$ in the disk, $L(h)>0$.

For an eigenfunction, $h\in H^{1/2}(\mathbb S^1)\cap C(\mathbb S^1)$
and $v-Ah\in H^1_0$. Its Dirichlet equation gives $v-Ah=Tv$,
so $v=SAh$. Finally,
$L(h)=\int qv=\mu\int_\Omega u=0$.
\end{proof}

\subsection{Multiplicity from boundary signs}

\begin{lemma}\label{lem:sa-multiplicity}
The real eigenspace at the first positive Neumann eigenvalue satisfies
\begin{equation}\label{eq:sa-multiplicity}
 \dim_{\mathbb R}E_\mu\le2.
\end{equation}
\end{lemma}
\begin{proof}
The boundary trace is injective on $E_\mu$: zero trace
would put the pullback in $H^1_0$, where
\eqref{eq:sa-form-bound} forces it to vanish. Its nonzero traces
take both signs, by the strict positivity of $L$ and $L(h)=0$.

Every nonzero eigenfunction in $E_\mu$ has exactly two nodal domains,
even when $\mu$ is multiple. To see this directly, if three nodal
domains existed, choose nonzero constants on the restrictions to
two of them so that their combination has zero mean. Each restriction,
extended by zero, lies in $H^1(\Omega)$ and has Rayleigh quotient
$\mu$. To justify the zero extension without boundary regularity, take a
positive nodal component $D$. The function
$(u-\epsilon)_+\mathbf1_D$, where $\mathbf1_D$ is the indicator of
$D$, is locally $H^1$: near every compact
piece of the interior boundary of $D$, continuity and $u=0$ make
it vanish for sufficiently small neighborhoods. Its global gradient
is $\mathbf1_{D\cap\{u>\epsilon\}}\nabla u$, and its function and
gradient norms are bounded by those of $u$. Letting $\epsilon\downarrow0$
in $H^1$ gives $u\mathbf1_D\in H^1$ with gradient
$\mathbf1_D\nabla u$. Negative components are treated with $-u$.
Testing the weak eigenfunction equation by this restriction gives
its Rayleigh quotient. The combination is therefore
a Neumann minimizer and an eigenfunction, but it vanishes on the third
open nodal domain. Interior analyticity of solutions of
$-\Delta u=\mu u$ forces it to vanish everywhere, a contradiction.
Both sign sets are consequently connected. The nodal restrictions
may equivalently be taken on the disk, where the pullback is continuous
and the energy and weighted mass identities are unchanged.

There cannot be four cyclically alternating strict signs in the
continuous disk boundary trace. Indeed, boundary continuity attaches
each such point to the corresponding interior sign set through a
small half-disk of the same sign. The connected positive and negative
sets contain paths joining the two positive and the two negative
points, respectively. After removing loops these can be taken as
simple arcs with interior in the disk. Two arcs joining alternating
boundary pairs intersect, by the Jordan separation theorem, whereas
the sign sets are disjoint. This argument takes place on the disk
and requires no regularity of $\partial\Omega$.

Parametrize the circle by $e^{i\theta}$. Suppose $\dim E_\mu\ge3$.
In its injective trace space choose a
nonzero $h$ such that
\[
 L(h\cos\theta)=L(h\sin\theta)=0.
\]
We already have $L(h)=0$. Cut the circle at a point of a negative
interval, and let $\theta_-,\theta_+$ be the infimum and supremum of
$\{\theta:h(\theta)>0\}$ in $(0,2\pi)$. Then
$0<\theta_-<\theta_+<2\pi$. Absence of four alternating signs implies
$h\ge0$ on $(\theta_-,\theta_+)$ and $h\le0$ on its complement.
The trigonometric polynomial
\[
 \varphi(\theta)=\cos\left(\theta-\frac{\theta_-+\theta_+}{2}\right)
              -\cos\left(\frac{\theta_+-\theta_-}{2}\right)
\]
has the same weak sign pattern and is strictly positive where
$h>0$. Thus $h\varphi$ is continuous, nonnegative, and nonzero, so
$L(h\varphi)>0$. But $\varphi$ belongs to the real span of $1,\cos\theta,
\sin\theta$, and all three moments vanish. This contradiction
proves \eqref{eq:sa-multiplicity}.
\end{proof}

\section{The weak divergence--curl form}
\label{sec:sa-hodge}

We retain the conformal hypothesis of Section~\ref{sec:sa-spectral}
and write $\mu=\mu_2(\Omega)$.
We use Rohleder's variational approach to Neumann gradients
\cite{Rohleder2023}. The decomposition below is formulated with weak
normal traces and therefore applies when the boundary has no classical
normal.

For a real vector field $X=(X_1,X_2)$ define
$\operatorname{curl}X=\partial_xX_2-\partial_yX_1$.
Let $\mathcal X(\Omega)$ consist of fields in $L^2(\Omega;\mathbb R^2)$
whose divergence and curl belong to $L^2$ and which have weak zero
normal trace, meaning
\begin{equation}\label{eq:sa-weak-normal}
 \int_\Omega X\cdot\nabla\eta
   +\int_\Omega (\operatorname{div}X)\eta=0
 \qquad(\eta\in H^1(\Omega)).
\end{equation}
For $X,Y\in\mathcal X(\Omega)$, set
\[
 \mathcal E_\mu(X,Y)=\int_\Omega
  \bigl(\operatorname{div}X\operatorname{div}Y
       +\operatorname{curl}X\operatorname{curl}Y-\mu X\cdot Y\bigr).
\]

\begin{lemma}\label{lem:sa-hodge}
Every $X\in\mathcal X(\Omega)$ has an
$L^2(\Omega;\mathbb R^2)$-orthogonal decomposition
\[
 X=\nabla f+\nabla^\perp h,\qquad
 \nabla^\perp h=(-h_y,h_x),
\]
where $f$ belongs to the Neumann operator domain, has mean zero,
and satisfies $\Delta f=\operatorname{div}X$, while
$h\in H^1_0(\Omega)$ belongs to the Dirichlet operator domain and
satisfies $\Delta h=\operatorname{curl}X$.
Write $\nabla E_\mu=\{\nabla u:u\in E_\mu\}$. Then
\begin{equation}\label{eq:sa-form-nonnegative}
 \mathcal E_\mu(X,X)\ge0,\qquad
 \mathcal E_\mu(X,X)=0
 \quad\Longleftrightarrow\quad X\in\nabla E_\mu.
\end{equation}
\end{lemma}
\begin{proof}
Taking $\eta=1$ in \eqref{eq:sa-weak-normal} gives the Neumann
compatibility condition. The Poincar\'e inequality on $\Omega$ and
Lax--Milgram therefore give $f$; the Dirichlet problem gives $h$.
The field $\nabla^\perp h$ has weak zero normal trace, because this
holds for compactly supported smooth approximants to $h$ and passes
to all $H^1$ tests by convergence of their gradients. The residual
\[
 Z=X-\nabla f-\nabla^\perp h
\]
has zero divergence, curl, and weak normal trace. Its components
are distributionally harmonic, hence smooth in the interior.
Simple connectivity gives a real smooth potential $P$ with
$\nabla P=Z$.

To use $P$ in \eqref{eq:sa-weak-normal}, fix an interior ball and
subtract the average of $P$ over that ball. Apply the Poincar\'e
inequality with this fixed average to
$P_M=\max\{-M,\min\{P,M\}\}$. For large $M$, the truncation agrees
with $P$ on the ball and has zero average there. Hence
$\|P_M\|_2\le C\|Z\|_2$. Fatou's lemma gives $P\in H^1(\Omega)$,
and \eqref{eq:sa-weak-normal}, tested with $P$, gives
$\int_\Omega|Z|^2=0$.

Testing the zero-normal-trace identity for $\nabla^\perp h$ with
$f$ proves orthogonality. Expansion in the Neumann and Dirichlet
eigenfunctions gives
\[
 \|\Delta f\|_2^2\ge\mu\|\nabla f\|_2^2,\qquad
 \|\Delta h\|_2^2\ge\lambda_1^D\|\nabla h\|_2^2.
\]
The strict inequality \eqref{eq:sa-gap} proves nonnegativity. Equality
forces $h=0$ and
$f\in E_\mu$, proving the nullspace statement.
\end{proof}

\subsection{The energy in conformal coordinates}

For $X=(X_1,X_2)$ and $Y=(Y_1,Y_2)$ in $\mathcal X(\Omega)$, put
$U=(X_1+iX_2)\circ\Phi$ and $V=(Y_1+iY_2)\circ\Phi$. Then
\[
 2U_z=F\bigl(\operatorname{div}X+i\operatorname{curl}X\bigr)\circ\Phi.
\]
Changing variables gives
\begin{equation}\label{eq:sa-conformal-form}
 \mathcal E_\mu(X,Y)=\mathfrak h_q(U,V)
 :=\operatorname{Re}\int_{\mathbb D}
       (4U_z\overline{V_z}-qU\overline V).
\end{equation}
Distributional differentiation and change of variables justify this
identity for fields in $\mathcal X(\Omega)$. For a real eigenfunction
$u\in E_\mu$, with $v=u\circ\Phi$,
\[
 W=(u_x+iu_y)\circ\Phi=\frac{2v_{\bar z}}{\overline F},\qquad
 W_z=-\frac\mu2Fv,\qquad
 \int_{\mathbb D}q|W|^2=\mu\int_\Omega|\nabla u|^2.
\]
The weak Neumann equation says that $\nabla u$ satisfies
\eqref{eq:sa-weak-normal}.

\begin{remark}[Boundary multiplication]
On the smooth domains $\Phi(\rho\mathbb D)$, $0<\rho<1$, a real
scalar multiple of a tangent boundary vector remains tangent. The
Dirichlet correction in the resolvent has zero trace, so solving
$(-\Delta-q_\rho)U=0$ on $\mathbb D$ with the multiplied data
preserves tangency. Here $q_\rho=\mu_\rho|\Phi_\rho'|^2$ is the
potential defined in the introduction. On the limiting domain, a
multiplied distributional trace need not define a field in
$\mathcal X(\Omega)$. Section~\ref{sec:allp2-spectral} projects the
smooth fields onto the gradients of all eigenfunctions whose
eigenvalues converge to $\mu$; their limits with vanishing energy
then lie in $\nabla E_\mu$.
\end{remark}

\section{Boundary kernels for integrable subcritical potentials}
\label{sec:holder-kernel-audit}

Throughout this section, let $q\ge0$ belong to $L^p(\mathbb D)$ for
some $1<p<\infty$. We call $q$ subcritical if, for some $0\le d<1$,
\begin{equation}\label{eq:hk-form}
 \int_{\mathbb D}q|w|^2\,dA
 \le d\int_{\mathbb D}|\nabla w|^2\,dA,
 \qquad w\in H^1_0(\mathbb D).
\end{equation}
Assume this condition throughout the section.
All constants are uniform when $p$, $d<1$, and an upper bound for
$\|q\|_p$ are fixed. Constants in H\"older estimates may also depend
on the chosen exponent; constants for a fixed interior point may
depend on its distance from $\mathbb S^1$.

Let $G=(-\Delta_{\mathbb D,D})^{-1}$, $Tf=G(qf)$, and $S=(I-T)^{-1}$
on $L^2(q)=L^2(\mathbb D,q\,dA)$. The bound $\|T\|\le d<1$, proved
below, justifies the inverse. For $x,y\in\mathbb D$, $x\ne y$, and
$s\in\mathbb S^1$, the ordinary disk Green and Poisson kernels are
\[
 \mathcal G_0(x,y)=\frac1{4\pi}
 \log\left(1+\frac{(1-|x|^2)(1-|y|^2)}{|x-y|^2}\right),
 \qquad
 P_s(x)=\frac{1-|x|^2}{2\pi|x-s|^2}.
\]
Boundary integrals use arclength measure $ds$. For $x\in\mathbb D$
and distinct $s,t\in\mathbb S^1$, define
\begin{align}
 K_q(x,s)&=P_s(x)+\sum_{j\ge1}(T^jP_s)(x),\label{eq:hk-K}\\
 R_q(s,t)&=\sum_{j\ge0}\int_{\mathbb D}
      P_s(z)(T^jP_t)(z)q(z)\,dA(z),\qquad
 N_q(s,t)=\frac1{\pi|s-t|^2}+R_q(s,t).
 \label{eq:hk-R}
\end{align}
Write $N_0(s,t)=1/(\pi|s-t|^2)$ for $s\ne t$.
Each term is an iterated integral of nonnegative kernels. This
interpretation does not require $P_s\in L^2(q)$.

\begin{theorem}[Kernel estimates]\label{thm:hk-estimates}
The kernels in \eqref{eq:hk-K}--\eqref{eq:hk-R} are finite for every
$x\in\mathbb D$, $s\in\mathbb S^1$, and every pair of distinct
$s,t\in\mathbb S^1$. They satisfy
\begin{align}
 P_s(x)&\le K_q(x,s)\le C P_s(x),\label{eq:hk-Pcomparison}\\
 0\le R_q(s,t)&\le C\|q\|_p|s-t|^{-2/p}\quad(s\ne t),
 \label{eq:hk-Rpoint}\\
 \iint_{\mathbb S^1\times\mathbb S^1}R_q(s,t)\,ds\,dt
 &\le\frac{\int_{\mathbb D}q\,dA}{1-d}.
 \label{eq:hk-Rmass}
\end{align}
The correction $R_q$ is symmetric and continuous off the diagonal.
For every fixed $z_0\in\mathbb D$,
\begin{equation}\label{eq:hk-kholder}
 K_q(z_0,\cdot)\in C^{0,\beta}(\mathbb S^1),\qquad
 0<\beta<\min\{1,2-2/p\}.
\end{equation}
The Dirichlet Green function $\mathcal G_q$ of $-\Delta-q$ satisfies
the following bounds for distinct $x,y\in\mathbb D$ and for $z$ in
a boundary collar separated from $z_0$:
\begin{equation}\label{eq:hk-Gbounds}
 0\le\mathcal G_q(x,y)\le C\log(1+|x-y|^{-1}),\qquad
 \mathcal G_q(z_0,z)\le C_{z_0}(1-|z|).
\end{equation}
\end{theorem}

\subsection{A positive solution and the Green function}

The variational equation for $Tf$, tested against $Tf$, gives
\[
 \|Tf\|_{L^2(q)}^2
 \le d\|\nabla Tf\|_2^2
 =d\langle f,Tf\rangle_{L^2(q)}.
\]
Thus $T$ is a positive self-adjoint contraction with $\|T\|\le d$;
positivity preservation follows also from $\mathcal G_0\ge0$.
Set $r=2p/(p+1)>1$. H\"older's inequality and the Dirichlet estimate
on the disk give
\begin{equation}\label{eq:hk-source}
 \|qf\|_r\le\|q\|_p^{1/2}\|f\|_{L^2(q)},\qquad
 T:L^2(q)\longrightarrow W^{2,r}\cap H^1_0
           \longrightarrow C(\overline{\mathbb D}).
\end{equation}
The resulting continuous functions vanish on $\mathbb S^1$.
For $q=0$ the stated kernels are the ordinary disk kernels. Assume
henceforth that $q\not\equiv0$ and define
\[
 \phi=S1=1+T\phi.
\]
The resolvent series and \eqref{eq:hk-source} imply
\begin{equation}\label{eq:hk-gauge-bounds}
 1\le\phi\le C,\qquad \phi|_{\mathbb S^1}=1,
 \qquad -\Delta\phi=q\phi,\qquad
 \phi-1\in W^{2,p}\cap H^1_0.
\end{equation}
The $W^{2,p}$ regularity follows by applying the Dirichlet estimate to
$q\phi\in L^p$. Let $I_2$ denote the $2\times2$ identity matrix and set
$\mathsf A=\phi^2I_2$. For every $0<\beta<\min\{1,2-2/p\}$,
\[
 \mathsf A\in C^{0,\beta}(\overline{\mathbb D}),
 \qquad I_2\le\mathsf A\le C I_2,\qquad \mathsf A|_{\mathbb S^1}=I_2.
\]
Also $\nabla\phi\in L^{p_0}$ for some $p_0>2$. Multiplication by $\phi$
and by $1/\phi$ is consequently bounded on $H^1(\mathbb D)$ and
preserves $H^1_0(\mathbb D)$: use
$H^1(\mathbb D)\hookrightarrow L^{2p_0/(p_0-2)}$ for
the product with $\nabla\phi$.

For smooth compactly supported $u,v$, the equation for $\phi$ gives
\begin{equation}\label{eq:hk-form-gauge}
 \int_{\mathbb D}\phi^2
       \nabla(u/\phi)\cdot\nabla(v/\phi)
 =\int_{\mathbb D}\nabla u\cdot\nabla v
       -\int_{\mathbb D}q uv.
\end{equation}
By density, the multiplication bounds and H\"older's inequality extend
the identity to all $u,v\in H^1_0(\mathbb D)$.
Let $\mathcal G_{\mathsf A}$ be the Dirichlet
Green kernel of $-\operatorname{div}(\mathsf A\nabla)$. Applying the
inverse operators in
\eqref{eq:hk-form-gauge}, first to smooth sources, gives
\begin{equation}\label{eq:hk-green-gauge}
 \mathcal G_q(x,y)=\phi(x)\phi(y)\mathcal G_{\mathsf A}(x,y).
\end{equation}
This identity initially holds almost everywhere. To obtain it at
every pair of distinct interior points, we identify continuous
representatives of the kernels.

For $x\in\mathbb D$, the section $\mathcal G_0(x,\cdot)$ belongs to
$L^2(q)$ and depends continuously on $x$ in that space. Indeed,
\[
 \mathcal G_0(x,z)=\frac1{2\pi}
 \bigl(\log|1-x\overline z|-\log|x-z|\bigr).
\]
On each compact set of poles, the first term is uniformly continuous
for $z\in\overline{\mathbb D}$. The logarithm belongs locally to
$L^{2p/(p-1)}$; multiplying it by a compactly supported cutoff equal
to one on the disk of radius two and using continuity of translations
gives continuity of the second term in $L^{2p/(p-1)}(\mathbb D)$.
H\"older's inequality then gives continuity
in $L^2(q)$. Since $S$ is bounded, the resolvent correction
\[
 \left\langle\mathcal G_0(x,\cdot),
 S\mathcal G_0(y,\cdot)\right\rangle_{L^2(q)}
\]
is jointly continuous for $x,y\in\mathbb D$, including $x=y$.
Its expansion into nonnegative Green chains follows from Tonelli's
theorem and $S=\sum_{j\ge0}T^j$; for fixed $x$, the correction series
converges uniformly in $y$ by \eqref{eq:hk-source} and $\|T\|\le d<1$.
Thus the chain representative of $\mathcal G_q$ is continuous off
the diagonal. The inverse-operator identity giving
\eqref{eq:hk-green-gauge} identifies it almost everywhere with the
continuous gauge representative, so the identity holds at every pair
of distinct interior points.

\subsection{Poisson kernel estimates}

Theorem 4.1, equations (27)--(28), of ter Elst--Ouhabaz
\cite{EO2019} applies to $\mathsf A$ on the disk with zero
zeroth-order potential. The coefficient class in that theorem requires
uniform ellipticity,
bounded coefficients, and a bounded $C^{0,\beta}$ seminorm; no
coefficient derivatives are required.
Their theorem gives, in dimension two,
\begin{equation}\label{eq:hk-EO}
 |\nabla_y\mathcal G_{\mathsf A}(x,y)|\le C|x-y|^{-1},
 \qquad
 |\nabla_x\nabla_y\mathcal G_{\mathsf A}(x,y)|\le C|x-y|^{-2},
\end{equation}
as well as a logarithmic bound for the Green function. These
derivatives extend locally $\beta$-H\"older continuously to
$\overline{\mathbb D}\times\overline{\mathbb D}$ off the diagonal.
The constants are uniform in the coefficient class just specified.
The boundary-collar bound in \eqref{eq:hk-Gbounds} follows by
integrating $\nabla_y\mathcal G_{\mathsf A}(z_0,y)$ from the nearest
boundary point and then using \eqref{eq:hk-green-gauge} and
\eqref{eq:hk-gauge-bounds}. The logarithmic bound follows directly
from \eqref{eq:hk-green-gauge}.

Let $\nu(s)$ be the outward unit normal to $\mathbb S^1$.
Proposition 5.5 of \cite{EO2019} identifies the Poisson kernel of
$-\operatorname{div}(\mathsf A\nabla)$ as
\[
 P_{\mathsf A}(x,s)
 =-\nu(s)\cdot\mathsf A(s)\nabla_y\mathcal G_{\mathsf A}(x,s)
 =-\partial_{\nu_s}\mathcal G_{\mathsf A}(x,s).
\]
It is nonnegative by the weak maximum principle. If
$\delta(x)\ge|x-s|/2$,
the first estimate in \eqref{eq:hk-EO} gives
$P_{\mathsf A}(x,s)\le C\delta(x)/|x-s|^2$. Otherwise, let $x_{\partial}=x/|x|$
be the nearest boundary point. The segment from $x_{\partial}$ to $x$ stays
at distance at least $|x-s|/2$ from $s$. Moreover
$\nabla_y\mathcal G_{\mathsf A}(x_{\partial},s)=0$: first differentiate
$\mathcal G_{\mathsf A}(x_{\partial},y)=0$ for interior $y$, then use the
continuous extension off the diagonal. The mixed derivative estimate
along this segment therefore gives the same bound. Hence
\begin{equation}\label{eq:hk-conductivity-P}
 0\le P_{\mathsf A}(x,s)
 \le C\frac{\delta(x)}{|x-s|^2}\le C P_s(x).
\end{equation}

For $f\in C^\infty(\mathbb S^1)$, let $Af$ be its ordinary Poisson
extension. Then $SAf$ has trace $f$ and satisfies
$(-\Delta-q)(SAf)=0$ weakly. The quotient $(SAf)/\phi$ satisfies
$-\operatorname{div}(\mathsf A\nabla((SAf)/\phi))=0$ with the same trace.
Consequently,
\begin{equation}\label{eq:hk-K-gauge}
 K_q(x,s)=\phi(x)P_{\mathsf A}(x,s)
\end{equation}
initially for almost every $s$. Reversing the iterated Green integrals
in \eqref{eq:hk-K} gives
\begin{equation}\label{eq:hk-K-Green}
 K_q(x,s)=P_s(x)+\int_{\mathbb D}
            \mathcal G_q(x,z)q(z)P_s(z)\,dA(z).
\end{equation}
In a boundary collar, \eqref{eq:hk-Gbounds} and
$\delta(z)P_s(z)\le1/\pi$ dominate the integrand by $C_xq(z)$,
uniformly in $s$. On the complementary compact disk, $P_s(z)$ is
uniformly bounded and $q(z)\mathcal G_q(x,z)$ is integrable by
H\"older's inequality and the logarithmic bound. Thus \eqref{eq:hk-K-Green} is
finite and continuous in $s$ by dominated convergence. The right
side of \eqref{eq:hk-K-gauge} is continuous as well, so the equality
holds at every boundary point.

Equations \eqref{eq:hk-conductivity-P}--\eqref{eq:hk-K-gauge} prove
the upper estimate in \eqref{eq:hk-Pcomparison}; the term $P_s$ in
\eqref{eq:hk-K} gives the lower estimate. Boundary H\"older continuity
of $\nabla_y\mathcal G_{\mathsf A}(z_0,s)$ proves
\eqref{eq:hk-kholder}.

\subsection{Bounds for the boundary correction}

Tonelli's theorem gives
\begin{equation}\label{eq:hk-R-integral}
 R_q(s,t)=\int_{\mathbb D}P_s(z)K_q(z,t)q(z)\,dA(z).
\end{equation}
Let $p'=p/(p-1)$. On the part of the disk defined by
$|z-s|\le|z-t|$, use
$\delta(z)\le|z-s|$ and
$|z-t|\ge\max\{|z-s|,|s-t|/2\}$ to obtain
$P_s(z)P_t(z)\le C(|s-t|+|z-s|)^{-2}$. Polar integration yields
\[
 \int_{\mathbb D}|P_sP_t|^{p'}\,dA
 \le C_p |s-t|^{2-2p'},
 \qquad \|P_sP_t\|_{p'}\le C_p |s-t|^{-2/p}.
\]
The complementary part has the same estimate. H\"older's inequality and
\eqref{eq:hk-Pcomparison} now give \eqref{eq:hk-Rpoint}.
For boundary pairs a positive distance from the diagonal,
$P_sP_t$ is uniformly bounded pointwise. Thus
\eqref{eq:hk-R-integral} and dominated convergence prove
off-diagonal continuity. Reversing the iterated Green integrals gives
symmetry.

The pointwise bound need not be integrable over $\mathbb S^1\times\mathbb S^1$ when
$p\le2$. To obtain the $L^1$ estimate, integrate the defining series:
\[
 \iint R_q(s,t)\,ds\,dt
 =\sum_{j\ge0}\langle1,T^j1\rangle_{L^2(q)}
 =\int_{\mathbb D}q\phi\,dA
 \le\frac{\int_{\mathbb D}q\,dA}{1-d}.
\]
This proves \eqref{eq:hk-Rmass} and completes the proof of
Theorem~\ref{thm:hk-estimates}.

\subsection{Regularity of the Hilbert-space curve}

Corollary~\ref{cor:rkchecked-extension} proves conditional negativity
of the quotient $K_q(z_0,s)K_q(z_0,t)/N_q(s,t)$, with zero diagonal.
Schoenberg's construction \cite{Schoenberg1938}, given in
Proposition~\ref{prop:rkchecked-realization}, represents this quotient
as squared distances in a real Hilbert space. The estimates above
give the regularity needed to multiply boundary traces.

\begin{proposition}\label{prop:hk-curve}
Fix $z_0\in\mathbb D$. Let $\mathcal H$ be a real Hilbert space and
suppose $b:\mathbb S^1\to\mathcal H$ satisfies
\[
 \|b(s)-b(t)\|^2=D_q(s,t)
 :=\frac{K_q(z_0,s)K_q(z_0,t)}{N_q(s,t)},\quad s\ne t,
 \qquad D_q(s,s)=0.
\]
Then $b$ is injective and belongs to $C^{1,\beta}(\mathbb S^1;\mathcal H)$
for every
\begin{equation}\label{eq:hk-gamma-range}
 0<\beta<\min\{\tfrac12,1-1/p\}.
\end{equation}
If $b(s_*)=0$ at one fixed $s_*\in\mathbb S^1$, its
$C^{1,\beta}$ norm is uniformly bounded.
In circle arclength coordinates,
\begin{equation}\label{eq:hk-speed}
 \left\|\frac{d}{d\theta}b(e^{i\theta})\right\|^2
 =\pi K_q(z_0,e^{i\theta})^2>0.
\end{equation}
\end{proposition}
\begin{proof}
Finiteness of $N_q$ and positivity of $K_q$ give injectivity, while
$N_q\ge1/(\pi|s-t|^2)$ gives the uniform Lipschitz bound.
Write $\kappa(\theta)=K_q(z_0,e^{i\theta})$.
The pointwise correction bound gives
\begin{align}
 D_q(e^{i\theta},e^{i\eta})
 &=\frac{\pi \kappa(\theta)\kappa(\eta)|e^{i\theta}-e^{i\eta}|^2}
 {1+\pi|e^{i\theta}-e^{i\eta}|^2R_q(e^{i\theta},e^{i\eta})}
 \notag\\
 &=\pi \kappa(\theta)\kappa(\eta)(\theta-\eta)^2
       +O(|\theta-\eta|^{4-2/p}).\label{eq:hk-metric-expansion}
\end{align}
The chord error is $O(|\theta-\eta|^4)$, which is included in the
displayed remainder since $p<\infty$.
Fix $\beta$ in the range \eqref{eq:hk-gamma-range}.
By \eqref{eq:hk-kholder}, $\kappa$ is $2\beta$-H\"older continuous.
Put $\widetilde b(\theta)=b(e^{i\theta})$. We use the identity
\[
 \|\widetilde b(\theta+h)-2\widetilde b(\theta)+\widetilde b(\theta-h)\|^2
 =2D_q(e^{i(\theta+h)},e^{i\theta})
  +2D_q(e^{i\theta},e^{i(\theta-h)})
  -D_q(e^{i(\theta+h)},e^{i(\theta-h)}).
\]
By \eqref{eq:hk-metric-expansion}, this gives
\begin{equation}\label{eq:hk-second-difference}
 \|\widetilde b(\theta+h)-2\widetilde b(\theta)+\widetilde b(\theta-h)\|
 \le C|h|^{1+\beta}.
\end{equation}
Indeed, replace the three products of $\kappa$-values by $\kappa(\theta)^2$,
at a cost $O(|h|^{2+2\beta})$. The quadratic terms then cancel
exactly, and the remaining metric errors have higher order.

We next deduce continuous differentiability from
\eqref{eq:hk-second-difference}. Periodic mollification gives smooth
Hilbert-valued functions $\widetilde b_\epsilon$ satisfying the same estimate. Define
$v_h^\epsilon(\theta)=
[\widetilde b_\epsilon(\theta+h)-\widetilde b_\epsilon(\theta)]/h$.
Equation \eqref{eq:hk-second-difference}, centered at $\theta+h/2$,
gives
\[
 \|v_h^\epsilon-v_{h/2}^\epsilon\|_\infty\le C|h|^\beta.
\]
Dyadic summation to the derivative of the smooth function yields
$\|(\widetilde b_\epsilon)'-v_h^\epsilon\|_\infty\le C_\beta|h|^\beta$.
For two mollification scales,
\[
 \|(\widetilde b_\epsilon)'-(\widetilde b_{\epsilon'})'\|_\infty
 \le2C_\beta|h|^\beta
       +\frac2{|h|}\|\widetilde b_\epsilon-\widetilde b_{\epsilon'}\|_\infty.
\]
Taking $h$ small first shows that the derivatives are uniformly
Cauchy as the mollification scales tend to zero. Their continuous
limit is $\widetilde b'$ by the fundamental theorem of calculus.
Writing $v_h(\theta)=[\widetilde b(\theta+h)-\widetilde b(\theta)]/h$, we obtain
$\|\widetilde b'-v_h\|_\infty\le C_\beta|h|^\beta$.
Also \eqref{eq:hk-second-difference} gives
$\|v_h(\theta+h)-v_h(\theta)\|\le C|h|^\beta$.
Combining these two inequalities proves that $\widetilde b'$ is
$\beta$-H\"older. Dividing \eqref{eq:hk-metric-expansion} by
$(\theta-\eta)^2$ and passing to the derivative limit proves
\eqref{eq:hk-speed}.
\end{proof}

\subsection{Convergence of the kernels}

\begin{proposition}\label{prop:hk-stability}
Suppose $q_n\ge0$, $q_n\to q$ strongly in $L^p(\mathbb D)$, and all
potentials satisfy \eqref{eq:hk-form} with the same $d<1$. For each
fixed $z_0\in\mathbb D$,
$K_{q_n}(z_0,\cdot)\to K_q(z_0,\cdot)$ uniformly. Moreover,
$N_{q_n}\to N_q$ uniformly on compact off-diagonal sets, and
$D_{q_n}\to D_q$ uniformly on $\mathbb S^1\times\mathbb S^1$.

Let $b_n,b$ be Hilbert-space realizations of $D_{q_n},D_q$ with
$b_n(s_*)=b(s_*)=0$ at the same point $s_*\in\mathbb S^1$. Then
\[
 c_{n,\xi}(s):=\langle b_n(s),b_n(\xi)\rangle
 \longrightarrow
 c_\xi(s):=\langle b(s),b(\xi)\rangle
 \quad\text{in }C^{1,\beta'}(\mathbb S^1)
\]
for each fixed $\xi\in\mathbb S^1$ and every
$0<\beta'<\beta<\min\{1/2,1-1/p\}$. No common ambient
Hilbert space is required.
\end{proposition}
\begin{proof}
The preceding estimates are uniform in $n$. The resolvent identity
for smooth boundary data, followed by Fubini, gives
\begin{equation}\label{eq:hk-K-difference}
 K_{q_n}(z_0,s)-K_q(z_0,s)
 =\int_{\mathbb D}\mathcal G_{q_n}(z_0,z)
       (q_n-q)(z)K_q(z,s)\,dA(z).
\end{equation}
This formula first identifies the boundary densities almost
everywhere. In the boundary collar the product
$\mathcal G_{q_n}(z_0,z)K_q(z,s)$ is uniformly bounded by
\eqref{eq:hk-Gbounds} and \eqref{eq:hk-Pcomparison}. On the
complementary compact disk it is bounded by a constant times
$1+|\log|z-z_0||$. It is therefore uniformly bounded in $L^{p'}$,
and dominated convergence makes the displayed integral continuous
in $s$. Thus the formula holds for every $s$, and
\[
 \sup_s|K_{q_n}(z_0,s)-K_q(z_0,s)|\le C_{z_0}\|q_n-q\|_p.
\]

The corresponding identity for the boundary kernels is
\begin{equation}\label{eq:hk-N-difference}
 N_{q_n}(s,t)-N_q(s,t)
 =\int_{\mathbb D}K_{q_n}(z,s)(q_n-q)(z)K_q(z,t)\,dA(z),
 \qquad s\ne t.
\end{equation}
Let $\Lambda_q$ be the weak Dirichlet-to-Neumann map for $-\Delta-q$.
For smooth boundary data $f,g$, integration against the harmonic
extension $Ag$ gives
\[
 \langle(\Lambda_q-\Lambda_0)f,g\rangle
 =-\int_{\mathbb D}q(SAf)(Ag)\,dA
 =-\iint_{\mathbb S^1\times\mathbb S^1}R_q(s,t)f(s)g(t)\,ds\,dt.
\]
The first identity uses $SAf-Af\in H^1_0$ and harmonicity of $Ag$.
This proves that the off-diagonal kernel of $\Lambda_q$ is $-N_q$.
Write $S_n=(I-Gq_n)^{-1}$. Subtracting the two
Dirichlet-to-Neumann forms gives $-\int(q_n-q)(S_nAf)(SAg)$.
Fubini's theorem applies to the kernel: its absolute integral over
the two boundary variables is bounded by
$\int|q_n-q|(S_n1)\phi<\infty$.
The kernels are continuous off the diagonal, so their almost
everywhere equality extends there pointwise.

The elementary bound $P_s(z)P_t(z)\le C|s-t|^{-2}$ gives
\[
 |N_{q_n}(s,t)-N_q(s,t)|
 \le C|s-t|^{-2}\|q_n-q\|_1.
\]
Together with $N_{q_n},N_q\ge1/(\pi|s-t|^2)$ and the preceding
bound for $K_{q_n}(z_0,\cdot)-K_q(z_0,\cdot)$, this gives
\[
 |D_{q_n}(s,t)-D_q(s,t)|
 \le C_{z_0}|s-t|^2\|q_n-q\|_p,
\]
including the zero diagonal by continuity.

Finally the Gram formula is
\[
 c_{n,\xi}(s)=\tfrac12\bigl(
 D_{q_n}(s,s_*)+D_{q_n}(\xi,s_*)-D_{q_n}(s,\xi)\bigr).
\]
It gives uniform convergence of $c_{n,\xi}$ to $c_\xi$.
Proposition~\ref{prop:hk-curve} gives uniform scalar $C^{1,\beta}$
bounds. Compact embedding into $C^{1,\beta'}$, together with
uniqueness of the uniform limit, proves the asserted convergence
for the full sequence.
\end{proof}

\section{Boundary traces of Neumann gradients}
\label{sec:primitive-trace-all-p2}

We represent the boundary trace of the harmonic remainder of a
Neumann gradient as $g_0+\partial_\theta g_1$, with
$g_0,g_1\in L^r(\mathbb S^1)$ for some $r>1$. This representation
allows multiplication by real $C^1$ functions and gives estimates
for Poisson extension on compact subdisks. The argument applies
whenever $F\in A^{2p}$ for some $p>1$.

\subsection{An integrability estimate for the reciprocal derivative}

We use the following consequence of the non-strict integral-means
estimate of Binder--Ivrii--Prause--Smirnov
\cite[Section 2.3, equation (27)]{BIPS2026}: for a bounded univalent
map with derivative $F$,
\begin{equation}\label{eq:allp2-inverse}
 I_\varepsilon(F):=\int_{\mathbb D}
       \delta^\varepsilon |F|^{-2}\,dA<\infty
       \qquad(\varepsilon>0).
\end{equation}
Apply their estimate with exponent $-2$ and auxiliary parameter
$1+\varepsilon/4$. The two associated parameters are $2-(-2)=4$
and $((1+\varepsilon/4)-1)4=\varepsilon$; the required inequality
is $1+\varepsilon/4>4/4$. Their estimate bounds the radial integral of
$(1-\rho)^{\varepsilon}\int_0^{2\pi}|F(\rho e^{i\theta})|^{-2}\,d\theta$
near $\rho=1$. Since $F$ has no zeros and $1/F$ is bounded on compact
subdisks, polar integration gives \eqref{eq:allp2-inverse}.
The argument below needs only
$I_\varepsilon(F)<\infty$ for one $0\le\varepsilon<1-1/p$;
\eqref{eq:allp2-inverse} supplies such an exponent for every $p>1$.

\subsection{An estimate for the analytic numerator}

Let $p>1$, let $F=\Phi'\in A^{2p}$, and let $v$ be the pullback
of a real Neumann eigenfunction with eigenvalue $\lambda$, so that
\[
 -\Delta v=\lambda |F|^2v,\qquad \partial_\nu v=0,
 \qquad v\in W^{2,p}(\mathbb D)\cap L^\infty.
\]
Set
\[
 W=\frac{2v_{\bar z}}{\overline F},\qquad
 T_\lambda=G(\lambda|F|^2\,\cdot),\qquad
 H=(I-T_\lambda)W.
\]
Since $-\Delta W=\lambda|F|^2W$, the function $H$ is harmonic.
Write $H=H_++\overline f$, where $H_+$ and $f$ are analytic and
$H_+(0)=0$. The estimates below are uniform when the eigenvalues,
$W^{2,p}$ norms, and $\|F\|_{2p}$ are bounded.
For the normalized eigenfunctions on $\Phi(\rho\mathbb D)$,
these bounds follow from Lemma~\ref{lem:sc-scalar-cluster} below.

\begin{lemma}[An estimate for the analytic numerator]
\label{lem:primitive-numerator}
The analytic function $h=Ff$ belongs to $A^{2p}$. Moreover,
$T_\lambda W$ and $H_+$ are uniformly bounded on the disk under the
uniform assumptions just stated. For triples $(F_n,v_n,\lambda_n)$
satisfying the same hypotheses, use a subscript $n$ for the corresponding
functions. If $F_n\to F$ strongly in
$L^{2p}$, $v_n\to v$ strongly in $W^{2,p}$, and
$\lambda_n\to\lambda$, then $T_{\lambda_n}W_n\to T_\lambda W$
strongly in $W^{2,p}$ and $H_{+,n}\to H_+$ strongly in
$W^{1,2p}$, hence uniformly on the closed disk. Also $h_n\to h$
strongly in $L^{2p}$.
\end{lemma}
\begin{proof}
For every $p>1$, disk Sobolev embedding gives
\[
 W^{2,p}(\mathbb D)\hookrightarrow W^{1,2p}(\mathbb D),
 \qquad \nabla v\in L^{2p}.
\]
For $1<p<2$, the usual first-derivative exponent is
$2p/(2-p)>2p$; for $p=2$ every finite exponent is available;
for $p>2$ the first derivatives are bounded.
Cancelling $\overline F$ gives
\[
 \lambda |F|^2W=2\lambda Fv_{\bar z}\in L^p.
\]
Hence the disk Dirichlet estimate gives
\[
 T_\lambda W\in W^{2,p}\cap H^1_0\subset C^0,
 \qquad \nabla T_\lambda W\in L^{2p}.
\]
Since
\[
 H_+'=W_z-(T_\lambda W)_z
       =-\frac\lambda2Fv-(T_\lambda W)_z,
\]
we have $H_+'\in A^{2p}$. The normalization $H_+(0)=0$ and
$2p>2$ give a uniform bound for $H_+$ by Morrey embedding.
Finally,
\begin{equation}\label{eq:primitive-actual-numerator}
 h=Ff=2v_z-F\overline{T_\lambda W+H_+}.
\end{equation}
The first term belongs to $L^{2p}$ and the second is $F$ times a
bounded function. This proves the membership and uniform bounds.

For the convergence assertion, the same embedding gives
$\nabla v_n\to\nabla v$ in $L^{2p}$, while
$v_n\to v$ uniformly. Thus
\[
 2\lambda_nF_n(v_n)_{\bar z}
 \longrightarrow2\lambda Fv_{\bar z}\quad\text{in }L^p,
\]
so the zero-Dirichlet corrections converge in $W^{2,p}$. The identity
for $H_+'$ then gives convergence of its derivatives in $L^{2p}$.
The analytic normalization $H_+(0)=0$ fixes the constant, so
$H_{+,n}\to H_+$ in $W^{1,2p}$ and uniformly. Finally
\eqref{eq:primitive-actual-numerator} gives $h_n\to h$ in $L^{2p}$.
\end{proof}

\subsection{A weighted integrability estimate}

Fix
\[
 r=\frac{2p}{p+1}\in(1,2),\qquad
 \varepsilon=\frac{p-1}{2p}>0,\qquad
 \sigma=\frac{\varepsilon r}{2}=\frac{r-1}{2}.
\]
The estimate \eqref{eq:allp2-inverse} is uniform for
$F_\rho(z)=\rho F(\rho z)$, because the circular means of $|1/F|^2$
are nondecreasing:
\[
 I_\varepsilon(F_\rho)\le\rho^{-2}I_\varepsilon(F)
 \le4I_\varepsilon(F),\qquad 1/2\le\rho<1.
\]

\begin{lemma}[A weighted $L^r$ estimate]
\label{lem:primitive-weighted-Lr}
With the exponents above,
\begin{equation}\label{eq:primitive-weighted-Lr}
 \left(\int_{\mathbb D}\delta^\sigma|f|^r\,dA\right)^{1/r}
 \le \|h\|_{L^{2p}} I_\varepsilon(F)^{1/2}.
\end{equation}
In particular, the weighted norm is uniformly bounded for the
harmonic remainders of the gradients in
Lemma~\ref{lem:sc-scalar-cluster}.
\end{lemma}
\begin{proof}
Since $f=h/F$ and $r=2p/(p+1)$,
\[
 \delta^\sigma|f|^r
 =\bigl(|h|^{2p}\bigr)^{1/(p+1)}
   \bigl(\delta^\varepsilon |F|^{-2}\bigr)^{p/(p+1)}.
\]
Apply H\"older with conjugate exponents $p+1$ and $(p+1)/p$,
then take the $r$th root.
\end{proof}

\subsection{Boundary values of a primitive}

\begin{lemma}[An $L^r$ boundary primitive]
\label{lem:primitive-boundary}
Let $1<r<\infty$, $0\le\sigma<r-1$, and let $f$ be analytic
in $\mathbb D$ with
\[
 M^r=\int_{\mathbb D}\delta^\sigma|f|^r\,dA<\infty.
\]
Let $Q$ be the analytic primitive $Q'=f$, $Q(0)=0$. Then there
exists $Q^*\in L^r(\mathbb S^1)$ such that
\begin{equation}\label{eq:primitive-radial-tail}
 \|Q^*-Q(\rho\,\cdot)\|_{L^r(\mathbb S^1)}
 \le C_{r,\sigma} M(1-\rho)^{(r-1-\sigma)/r},
 \qquad 1/2\le\rho<1.
\end{equation}
Also $\|Q^*\|_r\le C_{r,\sigma}M$.

If analytic functions $f_n$ converge to $f$ locally uniformly in
$\mathbb D$ and their weighted norms are uniformly bounded,
then the corresponding primitives satisfy
$Q_n^*\to Q^*$ strongly in $L^r(\mathbb S^1)$.
\end{lemma}
\begin{proof}
For $1/2\le\rho<R<1$, radial integration and Minkowski give
\[
 \|Q(R\,\cdot)-Q(\rho\,\cdot)\|_r
 \le\int_\rho^R\|f(t\,\cdot)\|_r\,dt.
\]
With $r'=r/(r-1)$, H\"older's inequality bounds the last expression by
\[
 \begin{aligned}
 &\left(\int_\rho^R(1-t)^\sigma
                       \|f(t\,\cdot)\|_r^r t\,dt\right)^{1/r}
 \left(\int_\rho^R(1-t)^{-\sigma/(r-1)}
                       t^{-1/(r-1)}\,dt\right)^{1/r'}\\
 &\hspace{20mm}\le C_{r,\sigma}M(1-\rho)^{(r-1-\sigma)/r}.
 \end{aligned}
\]
The second integral remains finite as $R\uparrow1$ exactly when
$\sigma<r-1$. This proves that the radial primitives
are Cauchy in $L^r$ and gives \eqref{eq:primitive-radial-tail}.
Analytic interior estimates bound $\|Q(\tfrac12\,\cdot)\|_r$ by
$C_{r,\sigma}M$, giving the bound for $\|Q^*\|_r$.

For the convergence assertion, fix $\rho<1$. Local uniform
convergence gives $Q_n(\rho\,\cdot)\to Q(\rho\,\cdot)$ uniformly.
The two radial tails are uniformly bounded by
$C(1-\rho)^{(r-1-\sigma)/r}$. Letting first $n\to\infty$ and then
$\rho\uparrow1$ proves strong $L^r$ convergence.
\end{proof}

For the exponents chosen above,
\[
 \frac{r-1-\sigma}{r}=\frac{p-1}{4p}>0.
\]
By Lemma~\ref{lem:primitive-numerator}, $H_+$ extends continuously
to $\overline{\mathbb D}$; we use the same notation for its boundary
values. Write
$W^{-1,r}(\mathbb S^1)=(W^{1,r'}(\mathbb S^1))^*$, where $r'=r/(r-1)$.
Differentiating the primitive, with $\zeta=e^{i\theta}$, gives
\[
 \partial_\theta\overline{Q(\rho\zeta)}
 =-i\rho\overline\zeta\,\overline{f(\rho\zeta)}.
\]
Taking $\rho\uparrow1$ in $W^{-1,r}$, we obtain the distributional
trace $g$ of $H=H_++\overline f$:
\begin{equation}\label{eq:primitive-full-trace}
 g=H_++i\zeta\partial_\theta\overline{Q^*}
       \in W^{-1,r}(\mathbb S^1;\mathbb C).
\end{equation}
Equivalently, it has the explicit representation
\begin{equation}\label{eq:primitive-trace-decomposition}
 g=g_0+\partial_\theta g_1,\qquad
 g_0=H_++\zeta\overline{Q^*},\qquad
 g_1=i\zeta\overline{Q^*},\qquad g_0,g_1\in L^r(\mathbb S^1).
\end{equation}
Indeed $\partial_\theta(i\zeta)=-\zeta$. Both $L^r$ norms are
uniformly bounded in the approximating family.
This also gives strong $W^{-1,r}$ convergence when the harmonic
remainders converge locally uniformly and satisfy the uniform bounds
above. Indeed,
Lemma~\ref{lem:primitive-boundary} gives strong $L^r$ convergence of
the primitives. The derivatives of the analytic parts $H_+$ are
uniformly bounded in $A^{2p}$, so those parts have uniformly H\"older
continuous boundary values; local convergence identifies their
uniform limit.

The harmonic remainder is recovered from its trace:
\[
 H=Ag.
\]
Indeed, $H(\rho z)=A[H(\rho\,\cdot)](z)$, and the radial traces
converge to $g$ in $W^{-1,r}$. The Poisson kernel is smooth in its
boundary variable on every compact subdisk, so passage to the limit
proves the identity. Poisson extension is injective on boundary
distributions: the $j$th Fourier coefficient at radius $\rho$ is
$\rho^{|j|}$ times the boundary coefficient.

\subsection{Multiplication of boundary distributions}
\label{sec:allp2-smoothing}

Use the representation $g=g_0+\partial_\theta g_1$ from
\eqref{eq:primitive-trace-decomposition}. All estimates below use
$\|g_0\|_{L^r}+\|g_1\|_{L^r}$ for these chosen representatives.
Neither representative need be tangent. Nonnegativity of the energy
is applied only to real multiples of the full tangent trace on the
smooth approximating domains.

\begin{lemma}[Multiplication and Poisson extension]
\label{lem:allp2-smoothing}
Let $\mathcal H$ be a real Hilbert space and
$c\in C^1(\mathbb S^1;\mathcal H)$. Then
\begin{equation}\label{eq:allp2-product}
 cg=(cg_0-(\partial_\theta c)g_1)+\partial_\theta(cg_1).
\end{equation}
Thus multiplication by $c$ is continuous from
$W^{-1,r}(\mathbb S^1;\mathbb C)$ to
$W^{-1,r}(\mathbb S^1;\mathcal H\otimes\mathbb R^2)$.
For every measurable $B\Subset\mathbb D$,
\begin{equation}\label{eq:allp2-local-poisson}
 \|A(cg)\|_{L^2(B;\mathcal H\otimes\mathbb R^2)}
 \le C_{B,r}\|c\|_{C^1(\mathcal H)}
       \bigl(\|g_0\|_{L^r}+\|g_1\|_{L^r}\bigr).
\end{equation}
The constant does not depend on the Hilbert dimension.
If $g_{j,n}\to g_j$ in $L^r$ for $j=0,1$,
$g_n=g_{0,n}+\partial_\theta g_{1,n}$, and $c_n\to c$ in $C^1$, then
$c_ng_n\to cg$ in $W^{-1,r}(\mathbb S^1;\mathcal H\otimes\mathbb R^2)$,
and their Poisson extensions converge strongly in $L^2$ on every
compact subdisk.
\end{lemma}
\begin{proof}
The distributional product rule proves \eqref{eq:allp2-product}.
Equivalently, for a Hilbert-valued test function $\psi$,
\[
 \|\langle c,\psi\rangle\|_{W^{1,r'}}
 \le C\|c\|_{C^1}\|\psi\|_{W^{1,r'}(\mathcal H)},
\]
which proves the dual-space statement. More directly,
\[
 A(cg)(z)=\int_{\mathbb S^1}P_s(z)c(s)g_0(s)\,ds
       -\int_{\mathbb S^1}\partial_\theta\bigl(P_s(z)c(s)\bigr)g_1(s)\,ds.
\]
The Poisson kernel and its first angular derivative are uniformly
bounded for $z\in B$. H\"older's inequality for Bochner integrals gives
\eqref{eq:allp2-local-poisson}, with a constant independent of
$\dim\mathcal H$. Applying the formula to a difference proves
convergence.
\end{proof}

\begin{remark}[Local inversion]
For a real scalar $c\in C^1$, the identity $cg=0$ implies
$\operatorname{supp}g\subset\{c=0\}$. On any open interval where
$c\ne0$, take a smooth compactly supported cutoff $\chi$.
Then $\chi/c$, extended by zero, is $C^1$, and associativity of
the dual multiplication operators gives
$\chi g=(\chi/c)(cg)=0$.
\end{remark}

\section{Smooth approximation of the first positive eigenspace}
\label{sec:allp2-spectral}

For $0<\rho<1$, let
\[
 \Phi_\rho(z)=\Phi(\rho z),\quad F_\rho(z)=\rho F(\rho z),\quad
 J_\rho=|F_\rho|^2,\quad \Omega_\rho=\Phi(\rho\mathbb D),\quad
 \mu_\rho=\mu_2(\Omega_\rho),\quad q_\rho=\mu_\rho J_\rho.
\]
Each $\Omega_\rho$ has analytic Jordan boundary. As in
Lemma~\ref{lem:sa-resolvent}, set $T_\rho=G(q_\rho\,\cdot)$ and
$S_\rho=(I-T_\rho)^{-1}$. The quadratic form
\eqref{eq:sa-conformal-form} becomes
\[
 \mathcal E_\rho(U)=\int_{\mathbb D}(4|U_z|^2-q_\rho|U|^2).
\]
Let $m=\dim E_\mu$. We approximate all $m$ eigenfunctions together,
since their eigenvalues may split on $\Omega_\rho$.
\subsection{Convergence of the scalar spectrum and eigenfunctions}

\begin{lemma}[Convergence of the first positive eigenspace]
\label{lem:sc-scalar-cluster}
As $\rho\uparrow1$, the weights satisfy
$J_\rho\to J=|F|^2$ strongly in $L^p$.
The scalar Neumann and Dirichlet eigenvalues on $\Omega_\rho$ converge,
with multiplicity, to their counterparts on $\Omega$.  In particular,
for $\rho$ sufficiently near one,
\begin{equation}\label{eq:sc-uniform-gaps}
 \mu_{m+2}(\Omega_\rho)-\mu_\rho\ge c_*,\qquad
 \lambda_1^D(\Omega_\rho)-\mu_\rho\ge c_*,\qquad
 \int q_\rho|w|^2\le d_*\int|\nabla w|^2\quad(w\in H^1_0),
\end{equation}
where $c_*>0$ and $d_*<1$ do not depend on $\rho$.

Along any sequence $\rho\uparrow1$, pass to a subsequence and choose real
eigenfunctions $u_{k,\rho}$, $1\le k\le m$, at eigenvalues
$\mu_{k+1}(\Omega_\rho)$, orthonormal in $L^2(\Omega_\rho)$, so that
$v_{k,\rho}=u_{k,\rho}\circ\Phi_\rho$ converge strongly in $H^1(\mathbb D)$
and $W^{2,p}(\mathbb D)$, uniformly on the closed disk, and smoothly
on compact subdisks, to $v_k=u_k\circ\Phi$, where
$u_1,\ldots,u_m$ form a real $L^2(\Omega)$-orthonormal basis of
$E_\mu$. The pullbacks have uniform $L^\infty(\mathbb D)$ and
$W^{2,p}(\mathbb D)$ bounds.
\end{lemma}

\begin{proof}
Analytic dilation converges strongly in every finite Bergman norm,
so $F_\rho\to F$ in $L^{2p}$ and $J_\rho\to J$ in $L^p$.
On the fixed Hilbert space $H^1(\mathbb D)$, let $\mathsf D$ and
$\mathsf M_\rho$ be the bounded self-adjoint operators representing
\[
 (v,w)\longmapsto\int\nabla v\cdot\nabla w,
 \qquad (v,w)\longmapsto\int J_\rho vw,
\]
respectively. H\"older's inequality and the compact embedding
$H^1(\mathbb D)\hookrightarrow L^{2p'}(\mathbb D)$ imply that
$\mathsf M_\rho$ is compact and converges in operator norm to the
operator defined by $J$. The operators $\mathsf D+\mathsf M_\rho$
are uniformly positive: Poincar\'e's inequality controls a function
up to a constant, and the uniform positive lower bound for $J_\rho$
on a fixed compact disk controls that constant. Thus the compact
positive operators
\[
 (\mathsf D+\mathsf M_\rho)^{-1/2}\mathsf M_\rho
 (\mathsf D+\mathsf M_\rho)^{-1/2}
\]
converge in operator norm. Their nonzero eigenvalues are
$(1+\mu_j(\Omega_\rho))^{-1}$.  The same argument on $H^1_0$ applies
to the Dirichlet problem. The min--max principle and convergence
of the spectral projections give convergence of the eigenvalues
and their associated subspaces. Since the limiting subspace is
finite dimensional, a subsequence of orthonormal eigenbases converges.
When $m>1$, this does not ensure that a prescribed limiting
eigenfunction is approximated by exact first eigenfunctions.

The fixed-disk equations are
\[
 -\Delta v_{k,\rho}=\mu_{k+1}(\Omega_\rho)J_{\rho}v_{k,\rho},
 \qquad \partial_\nu v_{k,\rho}=0.
\]
Uniform $H^1$ bounds give uniform $L^{p_1}$ bounds for every finite $p_1$.
Choose $p_1$ so large that the right side is bounded in some $L^{p_0}$,
$p_0>1$.  Disk Neumann regularity gives $W^{2,p_0}$ and thus a
uniform $L^\infty$ bound.  A second application gives uniform
$W^{2,p}$.  Compact embedding identifies uniform convergence
with the already known $H^1$ limit.  Consequently the right sides
converge strongly in $L^p$, and Neumann regularity, including
convergence of the disk means, gives strong $W^{2,p}$ convergence.
Interior smooth convergence follows from analytic convergence of
$F_\rho$ on compact subdisks and local elliptic regularity.
The limiting strict gaps give \eqref{eq:sc-uniform-gaps}.
\end{proof}

\subsection{Convergence of the gradient traces}

For $1\le k\le m$, set
\[
 \lambda_{k,\rho}=\mu_{k+1}(\Omega_\rho),\qquad
 W_{k,\rho}=\frac{2(v_{k,\rho})_{\bar z}}{\overline{F_\rho}},\qquad
 H_{k,\rho}=W_{k,\rho}-G(\lambda_{k,\rho}J_\rho W_{k,\rho}).
\]
The eigenvalue $\lambda_{k,\rho}$ in this correction makes
$H_{k,\rho}$ harmonic, even when $\lambda_{k,\rho}\ne\mu_\rho$.
For the limiting basis write
\[
 W_k=\frac{2(v_k)_{\bar z}}{\overline F},\qquad H_k=(I-T)W_k.
\]
Let $g_{k,\rho}$ and $g_k$ be the distributional boundary traces
of $H_{k,\rho}$ and $H_k$, respectively. Then
$H_{k,\rho}=Ag_{k,\rho}$ and $H_k=Ag_k$.
Lemmas~\ref{lem:primitive-numerator}
and~\ref{lem:primitive-boundary}, together with
\eqref{eq:allp2-inverse}, give representations
$g_k=g_{k,0}+\partial_\theta g_{k,1}$ and
\begin{equation}\label{eq:allp2-strong-traces}
 g_{k,\rho}=g_{k,\rho,0}+\partial_\theta g_{k,\rho,1},\qquad
 g_{k,\rho,j}\longrightarrow g_{k,j}\quad\hbox{in }L^{2p/(p+1)}(\mathbb S^1),
 \qquad j=0,1.
\end{equation}
Thus the full traces $g_{k,\rho}$, including their positive and
negative Fourier modes, converge to $g_k$ in
$W^{-1,2p/(p+1)}(\mathbb S^1)$.

\subsection{Convergence of the boundary kernels}

The potentials $q_\rho$ converge strongly to $q$ in $L^p$ and
satisfy the common strict form bound \eqref{eq:sc-uniform-gaps}.
Choose the Hilbert curves associated with $q_\rho$ and $q$ so that
$b_\rho(s_0)=b(s_0)=0$ at one fixed $s_0\in\mathbb S^1$.
Proposition~\ref{prop:hk-stability} gives uniform convergence of
the squared-distance kernels and, for every $\xi\in\mathbb S^1$,
$C^1$ convergence
\[
 c_{\rho,\xi}(s)=\langle b_\rho(s),b_\rho(\xi)\rangle
 \longrightarrow c_\xi(s)=\langle b(s),b(\xi)\rangle.
\]
The curves may take values in different Hilbert spaces, so we compare
their inner products. The potentials need not be monotone.

\subsection{A point-evaluation estimate from the harmonic remainder}

Consider the gradient fields associated with the $m$ Neumann
eigenvalues converging to $\mu$. Their $L^2(\Omega_\rho)$
normalizations and the corresponding remainders are
\[
 e_{k,\rho}=\frac{\nabla u_{k,\rho}}{\sqrt{\lambda_{k,\rho}}},
 \qquad L_{k,\rho}=(I-T_\rho)(e_{k,\rho}\circ\Phi_\rho).
\]
The $L_{k,\rho}$ need not be harmonic when
$\lambda_{k,\rho}\ne\mu_\rho$. They converge locally uniformly to
\[
 L_k=\frac{H_k}{\sqrt\mu},
\]
which are harmonic and linearly independent over the reals.

To prove this convergence, observe that
\[
 q_\rho W_{k,\rho}
 =2\mu_\rho F_\rho(v_{k,\rho})_{\bar z}
 \longrightarrow 2\mu F(v_k)_{\bar z}=qW_k
 \quad\hbox{in }L^{2p/(p+1)}(\mathbb D).
\]
Here $F_\rho\to F$ in $L^{2p}$ and $v_{k,\rho}\to v_k$ in $H^1$.
The disk Dirichlet estimate gives $T_\rho W_{k,\rho}\to TW_k$ in
$W^{2,2p/(p+1)}$ and uniformly. Interior convergence of the gradients
and $\lambda_{k,\rho}\to\mu$ then give the stated local convergence
of the normalized remainders $L_{k,\rho}$.

\begin{lemma}[Uniform norm and pointwise estimates]
\label{lem:sc-point-evaluation}
Fix $z_0\in\mathbb D$ and a nonempty disk
$B\Subset\mathbb D$.  There is a constant $C$ independent of
$\rho$ near one and of the real Hilbert coefficient space such that
every field $X$ in the form domain of Section~\ref{sec:sa-hodge}
on $\Omega_\rho$, with vanishing normal trace and pullback
$U=X\circ\Phi_\rho$, satisfies
\begin{equation}\label{eq:sc-mass-coercivity}
 \|X\|_{L^2(\Omega_\rho)}^2
 \le C\bigl(\|(I-T_\rho)U\|_{L^2(B)}^2+
                        \mathcal E_\rho(U)\bigr).
\end{equation}
If in addition $(-\Delta-q_\rho)U=0$ componentwise, then
\begin{equation}\label{eq:sc-evaluation-coercivity}
 \|U(z_0)\|^2
 \le C\bigl(\|(I-T_\rho)U\|_{L^2(B)}^2+
                        \mathcal E_\rho(U)\bigr).
\end{equation}
\end{lemma}

\begin{proof}
The Hodge spectral decomposition gives
\[
 X=\sum_{k=1}^m a_k e_{k,\rho}+Y,
 \qquad Y\perp e_{k,\rho}\quad(1\le k\le m).
\]
For Hilbert-valued fields the coefficients $a_k$ are Hilbert
vectors.  The Neumann and Dirichlet spectral gaps in
\eqref{eq:sc-uniform-gaps}, together with the nonnegative energy
of the $m$ gradient fields, imply
\begin{equation}\label{eq:sc-residual-gap}
 \|Y\|_{L^2(\Omega_\rho)}^2\le C\mathcal E_\rho(U).
\end{equation}
Indeed, the eigenvalues of the form on gradient fields are
$\mu_j(\Omega_\rho)-\mu_\rho$, and on rotated Dirichlet gradients
are $\lambda_j^D(\Omega_\rho)-\mu_\rho$.

If $\widehat Y=Y\circ\Phi_\rho$, then
$\|\widehat Y\|_{L^2(q_\rho)}^2=\mu_\rho\|Y\|_2^2$.
On $B$, $q_\rho$ has a uniform positive lower bound, while
\[
 \|T_\rho\widehat Y\|_{C^0}
 \le C\|\widehat Y\|_{L^2(q_\rho)}.
\]
The Hilbert-valued version follows directly from the positive disk
Green kernel and Cauchy--Schwarz:
\[
 \|T_\rho\widehat Y(z)\|\le
 \left(\int \mathcal G_0(z,w)^2q_\rho(w)\,dA(w)\right)^{1/2}
 \|\widehat Y\|_{L^2(q_\rho)}.
\]
The first factor is bounded uniformly, since the logarithmic Green
singularity has uniformly bounded $L^{2p'}$ norm and $q_\rho$ is
bounded in $L^p$. Together with \eqref{eq:sc-residual-gap}, this bounds
$\|(I-T_\rho)\widehat Y\|_{L^2(B)}$ by
$C\sqrt{\mathcal E_\rho(U)}$.

The limiting $L_k$ are linearly independent on every nonempty $B$.
If a real combination were zero there, harmonic unique continuation
would make the corresponding $H$ zero throughout the disk.  The
associated gradient would satisfy $W=TW$ in $L^2(q)$,
contradicting $\|T\|<1$ unless the eigenfunction were zero.
Therefore the real Gram matrix
\[
 \left(\operatorname{Re}\int_B L_k\overline{L_\ell}\right)_{k,\ell=1}^m
\]
is positive definite. The smallest eigenvalues of the corresponding
Gram matrices for $L_{k,\rho}$ have a uniform positive lower bound
for $\rho$ near one. Otherwise a violating sequence would have a subsequence
as in Lemma~\ref{lem:sc-scalar-cluster}, with a positive definite
limiting Gram matrix, a contradiction.
From
\[
 (I-T_\rho)U=\sum_k a_kL_{k,\rho}+(I-T_\rho)\widehat Y
\]
we obtain
\[
 \sum_k\|a_k\|^2
 \le C\bigl(\|(I-T_\rho)U\|_{L^2(B)}^2+
                  \mathcal E_\rho(U)\bigr).
\]
Orthogonality and \eqref{eq:sc-residual-gap} prove
\eqref{eq:sc-mass-coercivity}. The constants do not depend on the
Hilbert-space dimension.

For the last assertion use an interior disk containing $z_0$.
The potentials $q_\rho$ and all their derivatives are uniformly
bounded there.  The interior elliptic point estimate for a
Helmholtz solution bounds $\|U(z_0)\|$ by its local $L^2$ norm.
The positive lower bound for $J_\rho$ on this disk bounds that local
norm by $\|X\|_{L^2(\Omega_\rho)}$. Summing the scalar squared
estimates over an orthonormal basis
gives the vector-valued estimate with the same constant.  This proves
\eqref{eq:sc-evaluation-coercivity}.
\end{proof}

\subsection{Approximation of a prescribed eigenfunction}
\label{sec:sc-corrected-approximation}

Fix $u=\sum_{k=1}^ma_k u_k\in E_\mu$, and let
$W=(u_x+iu_y)\circ\Phi$. Define
\[
 W_\rho^{\mathrm{cl}}=\sum_ka_kW_{k,\rho},\qquad
 g_\rho=\sum_ka_kg_{k,\rho},\qquad
 \widehat W_\rho=S_{\rho}Ag_\rho.
\]
The field $\widehat W_\rho$ solves $(-\Delta-q_\rho)\widehat W_\rho=0$.
Its difference $Z_\rho=\widehat W_\rho-W_\rho^{\mathrm{cl}}$ solves
\begin{equation}\label{eq:allp2-repair-equation}
 \begin{cases}
 (-\Delta-q_\rho)Z_\rho=-\displaystyle\sum_ka_k
       (\lambda_{k,\rho}-\mu_\rho)J_{\rho}W_{k,\rho},&\text{in }\mathbb D,\\
 Z_\rho=0,&\text{on }\mathbb S^1.
 \end{cases}
\end{equation}
With $r=2p/(p+1)>1$, weighted H\"older gives
\[
 \|J_\rho W_{k,\rho}\|_r
 \le \|F_\rho\|_{2p}\|F_\rho W_{k,\rho}\|_2
 =\|F_\rho\|_{2p}\sqrt{\lambda_{k,\rho}}.
\]
Since $\lambda_{k,\rho}-\mu_\rho\to0$, the right side of
\eqref{eq:allp2-repair-equation} tends to zero in $L^r$. The uniform
Dirichlet gap and disk regularity give $Z_\rho\to0$ in
$H^1_0(\mathbb D)$ and $W^{2,r}(\mathbb D)$, and uniformly on
$\overline{\mathbb D}$. Indeed, after the $H^1_0$
estimate the remaining source term satisfies
$\|q_\rho Z_\rho\|_{2p/(p+1)}
\le\|q_\rho\|_p\|Z_\rho\|_{2p/(p-1)}=o(1)$,
which gives the asserted Dirichlet regularity. Since $Z_\rho=0$
on $\mathbb S^1$, the field $\widehat W_\rho$ has the same boundary
values as $W_\rho^{\mathrm{cl}}$ and is therefore tangent.

The correction also tends to zero in the norm of the quadratic form, since
\[
 4\|(Z_\rho)_z\|_2^2+\int_{\mathbb D}q_\rho|Z_\rho|^2
 \le4\|\nabla Z_\rho\|_2^2+
       \|q_\rho\|_p\|Z_\rho\|_{2p'}^2\longrightarrow0.
\]
The field $Z_\rho\circ\Phi_\rho^{-1}$ has zero boundary value
and belongs to the form domain with vanishing normal trace.
Orthogonality of the scalar eigenfunctions gives
\[
 \mathcal E_\rho(W_\rho^{\mathrm{cl}})
 =\sum_ka_k^2\lambda_{k,\rho}(\lambda_{k,\rho}-\mu_\rho)\longrightarrow0.
\]
Nonnegativity and Cauchy--Schwarz for the form then imply, for
every fixed $z_0\in\mathbb D$,
\begin{equation}\label{eq:allp2-repair-energy}
 e_\rho:=\mathcal E_\rho(\widehat W_\rho)\longrightarrow0,
 \qquad \widehat W_\rho(z_0)\longrightarrow W(z_0).
\end{equation}
The two $L^r$ representatives of $g_\rho$ converge strongly by
\eqref{eq:allp2-strong-traces}. This approximates every $u\in E_\mu$,
even when the eigenvalues $\lambda_{k,\rho}$ are distinct. The
correction is needed because $(I-T_\rho)W_\rho^{\mathrm{cl}}$ need
not be harmonic.

\section{The Gaussian multiplier estimate}
\label{sec:coherent-near-null}

We first work on one smooth domain $\Omega_\rho$ and suppress
$\rho$ from the notation. For $f\in H^{1/2}(\mathbb S^1;\mathbb C)$, set
\[
 V_f=SAf,\qquad
 \mathcal E(U,V)=\mathfrak h_q(U,V),\qquad
 \mathcal E(U)=\mathfrak h_q(U,U),
\]
where $\mathfrak h_q$ is defined in \eqref{eq:sa-conformal-form}.
If $f$ is the trace of a tangent vector field, then
$\mathcal E(V_f)\ge0$. Multiplication by a real function preserves
this tangency. For data in the complexification of a real Hilbert
space, apply
$SA$ componentwise and sum the energies over an orthonormal basis.
The resulting form is independent of the basis, and all estimates
below are independent of the Hilbert-space dimension.

Fix $z_0\in\mathbb D$. Let $b:\mathbb S^1\to\mathcal H$ be the embedding
in Proposition~\ref{prop:hk-curve}, with $b(s_0)=0$ at a fixed
$s_0\in\mathbb S^1$. Write $K=K_q$, $N=N_q$, and
\[
 D_q(s,t)=\|b(s)-b(t)\|^2,\qquad
 N(s,t)D_q(s,t)=K(z_0,s)K(z_0,t).
\]
\subsection{The scalar commutator and its normalization}
\label{sec:sa-commutator}

Recall that $A$ is Poisson extension with respect to $d\theta$.
For $f,g\in H^{1/2}(\mathbb S^1;\mathbb C)$, define
\[
 \mathcal F_q(V_f,V_g)=\operatorname{Re}\int_{\mathbb D}
       (\nabla V_f\cdot\overline{\nabla V_g}-qV_f\overline{V_g}).
\]
With $R_q$ as in \eqref{eq:hk-R}, the resolvent expansion gives,
for bounded boundary data,
\begin{equation}\label{eq:sa-R-pairing}
 \langle Af,SAg\rangle_{L^2(q)}
   =\iint_{\mathbb S^1\times\mathbb S^1}R_q(s,t)f(s)\overline{g(t)}\,ds\,dt.
\end{equation}
For nonnegative data the formula follows term by term from Tonelli;
general bounded data follow by absolute domination. Indeed,
\[
 \iint R_q(s,t)\,ds\,dt=\langle1,S1\rangle_{L^2(q)}<\infty.
\]
Self-adjointness makes $R_q$ symmetric, and $N=N_0+R_q$.

\begin{lemma}\label{lem:sa-commutator}
For a real $c\in C^1(\mathbb S^1)$ and smooth complex $f$,
\begin{equation}\label{eq:sa-commutator}
 \begin{split}
 &\mathfrak h_q(V_{cf},V_{cf})-\mathfrak h_q(V_f,V_{c^2f})\\
 &\qquad=\frac12\iint_{\mathbb S^1\times\mathbb S^1}
 N(s,t)(c(s)-c(t))^2
       \operatorname{Re}\bigl(f(s)\overline{f(t)}\bigr)\,ds\,dt.
 \end{split}
\end{equation}
\end{lemma}
\begin{proof}
Testing the Helmholtz equation for $V_f$ against $V_g-Ag\in H^1_0$,
and using harmonic orthogonality, gives
\[
 \mathcal F_q(V_f,V_g)
 =\operatorname{Re}\int\nabla Af\cdot\overline{\nabla Ag}
       -\operatorname{Re}\langle Af,SAg\rangle_{L^2(q)}.
\]
The harmonic energy equals
\[
 \frac12\operatorname{Re}\iint N_0(s,t)
       (f(s)-f(t))\overline{(g(s)-g(t))}\,ds\,dt.
\]
The normalization may be checked on $f(e^{i\theta})=e^{i\theta}$:
both expressions equal $2\pi$. Combining this with
\eqref{eq:sa-R-pairing} and expanding the two quadratic terms
proves \eqref{eq:sa-commutator} with $\mathcal F_q$ in place of
$\mathfrak h_q$.

For $U=X+iY$ one has the pointwise identity
\[
 4|U_z|^2=|\nabla U|^2+2(X_xY_y-X_yY_x).
\]
Green's formula and polarization give
\[
 \mathfrak h_q(V_f,V_g)=\mathcal F_q(V_f,V_g)
       +\operatorname{Im}\int_{\mathbb S^1}\bar f\,\partial_\theta g.
\]
The boundary term is real symmetric by integration by parts. The
formula extends to $H^1$ lifts by smooth approximation and the
$H^{1/2}$ boundary pairing. For real $c$, its values at $(cf,cf)$
and $(f,c^2f)$ both equal
\[
 \int_{\mathbb S^1}c^2\operatorname{Im}(\bar f\,\partial_\theta f).
\]
The boundary terms therefore cancel, proving
\eqref{eq:sa-commutator}.
\end{proof}

For a $C^1$ multiplier with values in a real Hilbert space, sum
the identity over its coordinates. Finite orthogonal projections
converge uniformly on the compact ranges of the multiplier and its
derivative, hence in $C^1$. On each fixed smooth approximation the
corresponding lifts converge in $H^1$. On the kernel side, the squared
Lipschitz difference cancels the singularity of $N_0$, and the
remaining part is dominated by the integrable kernel $R_q$.
The resulting identity therefore holds in the Hilbert space,
with no constant depending on its dimension.

\subsection{The normalized exponential curve}

Consider the real symmetric Fock space
$\bigoplus_{n\ge0}\mathcal H^{\otimes_s n}$, where
$\mathcal H^{\otimes_s n}$ is the $n$th symmetric Hilbert tensor power of
$\mathcal H$ and $\mathcal H^{\otimes_s0}=\mathbb R$.
We use $x^{\otimes0}=1$ and $\|x^{\otimes n}\|=\|x\|^n$.
Write $e_0$ for the degree-zero unit vector. For $x\in\mathcal H$
and $t\ge0$, define
\[
 C_t(x)=e^{-t\|x\|^2/2}
       \bigoplus_{n\ge0}\frac{t^{n/2}}{\sqrt{n!}}x^{\otimes n},
 \qquad \mathcal B_t=C_t\circ b,\qquad
 d_t=e^{-t\|b\|^2/2}.
\]
Direct summation gives, for $x,y\in\mathcal H$,
\begin{equation}\label{eq:coherent-Gram}
 \|C_t(x)\|=1,\qquad
 \langle C_t(x),C_t(y)\rangle=e^{-t\|x-y\|^2/2}.
\end{equation}

For Hilbert-valued functions on $\mathbb S^1$, use the norm
$\|f\|_{C^1}=\|f\|_\infty+\|\partial_\theta f\|_\infty$.

\begin{lemma}[Uniform $C^1$ estimates]\label{lem:coherent-C1gamma}
For $b\in C^1(\mathbb S^1;\mathcal H)$, there is a constant $C_b$
depending only on $\|b\|_{C^1}$, not on the Hilbert dimension, such
that for $0\le t\le1$,
\[
 \|\mathcal B_t-e_0\|_{C^1}\le\sqrt t\,\|b\|_{C^1},
 \qquad \|d_tb-b\|_{C^1}\le C_b t.
\]
\end{lemma}
\begin{proof}
The defining series and its first derivatives converge uniformly on
bounded balls. Since $C_t(0)=e_0$, \eqref{eq:coherent-Gram} gives
\[
 \|C_t(x)-e_0\|^2=2(1-e^{-t\|x\|^2/2})\le t\|x\|^2,
 \qquad \|DC_t(x)h\|^2=t\|h\|^2.
\]
The chain rule proves the first estimate. For the second, use
$|d_t-1|\le t\|b\|^2/2$ and
$\partial_\theta d_t=-t d_t\langle b,\partial_\theta b\rangle$.
\end{proof}

\subsection{The Gaussian multiplier identity}

\begin{proposition}[Gaussian multiplier identity]
\label{prop:coherent-identity}
For every smooth tangent trace $f$ and every $t\ge0$,
\begin{equation}\label{eq:coherent-identity}
 \mathcal E(V_{\mathcal B_t f},V_{\mathcal B_t f})-\mathcal E(V_f,V_f)
   =\frac12\int_0^t\|V_{\mathcal B_\tau f}(z_0)\|^2
             \,d\tau.
\end{equation}
In particular, the left side is nonnegative.
\end{proposition}

\begin{proof}
Sum \eqref{eq:sa-commutator} over the coordinates of
\(\mathcal B_t\).  Since $\|\mathcal B_t(s)\|^2=1$, the mixed term equals
\(\mathcal E(V_f,V_f)\).  The remaining term is
\[
 \frac12\iint N(s,t')\|\mathcal B_t(s)-\mathcal B_t(t')\|^2
       \operatorname{Re}(f(s)\overline{f(t')})\,ds\,dt'.
\]
By \eqref{eq:coherent-Gram},
\[
 \frac12N(s,t')\|\mathcal B_t(s)-\mathcal B_t(t')\|^2
 =N(s,t')\bigl(1-e^{-tD_q(s,t')/2}\bigr)
 =\frac12K(z_0,s)K(z_0,t')\int_0^t e^{-\tau D_q(s,t')/2}\,d\tau.
\]
The inner product in \eqref{eq:coherent-Gram} then identifies the
integral with \eqref{eq:coherent-identity}.

The Hilbert-valued extension of Lemma~\ref{lem:sa-commutator}
justifies the coordinate sum. Its bounds are uniform for $\tau$ in
a bounded interval: the squared Lipschitz difference cancels the
singularity of $N_0$, while $R_q=N-N_0$ is integrable. They also
justify the integration in $\tau$.
\end{proof}

\subsection{The estimate with an \texorpdfstring{$L^r$}{Lr} primitive}

\begin{theorem}[Gaussian multiplier estimate]
\label{thm:allp2-gaussian}
On the smooth approximations of Section~\ref{sec:allp2-spectral},
fix $r>1$. There exist $\rho_0\in(0,1)$ and $C>0$, depending
only on $\Phi$, $p$, $z_0$, and $r$, with the following property.
For $\rho\in(\rho_0,1)$, let $g$ be a smooth tangent trace satisfying
$g=g_0+\partial_\theta g_1$ in the sense of distributions, where
$g_0,g_1\in L^r(\mathbb S^1;\mathbb C)$. Set
\[
 e=\mathcal E(V_g),\qquad M=\|g_0\|_r+\|g_1\|_r.
\]
Then
\begin{equation}\label{eq:allp2-gaussian}
 \mathcal E(V_{d_tbg})
 \le C\left(\frac et+|V_g(z_0)|^2+tM^2\right),\qquad 0<t\le1.
\end{equation}
The constant $C$ is independent of $\rho$, $g$, $t$, and the
Hilbert dimension.
\end{theorem}
\begin{proof}
For the curve $\mathcal B_t$ defined above, put
\[
 y(t)=\mathcal E(V_{\mathcal B_tg}),\qquad Z_t=V_{(\mathcal B_t-e_0)g}.
\]
Proposition~\ref{prop:coherent-identity} gives $y(0)=e$ and
$y'(t)=\tfrac12\|V_{\mathcal B_tg}(z_0)\|^2$.
Fix a nonempty disk $B\Subset\mathbb D$.
Lemma~\ref{lem:coherent-C1gamma} gives
$\|\mathcal B_t-e_0\|_{C^1}\le C\sqrt t$, so
Lemma~\ref{lem:allp2-smoothing} yields
\[
 \|(I-T)Z_t\|_{L^2(B)}
 =\|A((\mathcal B_t-e_0)g)\|_{L^2(B)}\le C\sqrt t\,M.
\]
The pointwise estimate \eqref{eq:sc-evaluation-coercivity}
therefore gives
\[
 \|Z_t(z_0)\|^2\le C\bigl(\mathcal E(Z_t)+tM^2\bigr).
\]
Nonnegativity of the form implies $\mathcal E(Z_t)\le2y(t)+2e$.
As $V_{\mathcal B_tg}(z_0)=e_0V_g(z_0)+Z_t(z_0)$, we obtain
\[
 y'(t)\le C\bigl(|V_g(z_0)|^2+y(t)+e+tM^2\bigr).
\]
Gronwall's inequality on $[0,1]$ yields
$y(t)\le C(e+t|V_g(z_0)|^2+t^2M^2)$.
The degree-one tensor component of $\mathcal B_tg$ is $\sqrt t\,d_tbg$.
All coordinate fields are tangent and their energies are nonnegative,
so $t\mathcal E(V_{d_tbg})\le y(t)$, proving the assertion.
\end{proof}

\paragraph{The exponent range.}
For the traces in Section~\ref{sec:primitive-trace-all-p2}, take
$r=2p/(p+1)>1$. The radial convergence exponent in
Lemma~\ref{lem:primitive-boundary} is $(p-1)/(4p)>0$, and
Proposition~\ref{prop:hk-curve} bounds $b_\rho$ uniformly in $C^1$.
Thus the multiplier estimate applies for every $p>1$.

\subsection{Limits of the multiplied fields}
\label{sec:scalar-multiplier-limits}

On the limiting domain $\Omega$, use $T=G(\mu|F|^2\,\cdot)$.
Let $u\in E_\mu$ be nonzero, set $W=(u_x+iu_y)\circ\Phi$,
and let $g=g_u$ be the boundary trace of $(I-T)W$ from
\eqref{eq:primitive-full-trace}.
Suppose that $W(z_0)=0$, and set $r=2p/(p+1)$. Choose
$\rho_n\uparrow1$ and use the subscript $n$ for the corresponding
quantities. The fields constructed in
Section~\ref{sec:sc-corrected-approximation} satisfy $e_n\to0$ and
$\widehat W_n(z_0)\to0$, and their two $L^r$ trace representatives
are uniformly bounded. Set
\[
 t_n=\max\{\sqrt{e_n},1/n\}.
\]
Then $t_n\to0$ and $e_n/t_n\le\sqrt{e_n}$, including when $e_n=0$.
After discarding finitely many terms, $t_n\le1$. With $V_h=S_nAh$
on the $n$th approximation, Theorem~\ref{thm:allp2-gaussian} gives
\begin{equation}\label{eq:allp2-diagonal}
 \mathcal E_n\bigl(V_{e^{-t_n\|b_n\|^2/2}b_ng_n}\bigr)\longrightarrow0.
\end{equation}

Choose a countable dense set of $\xi\in\mathbb S^1$. By
Proposition~\ref{prop:hk-stability}, the functions
\[
 c_{n,\xi}(s)=\langle b_n(s),b_n(\xi)\rangle
 \longrightarrow c_\xi(s)=\langle b(s),b(\xi)\rangle
 \quad\hbox{in }C^1.
\]
The vectors $b_n(\xi)$ have uniformly bounded norm. The extension
$V_h=S_nAh$ commutes with constant Hilbert-space projections, so
nonnegativity of the quadratic form gives
\[
 \mathcal E_n\bigl(V_{e^{-t_n\|b_n\|^2/2}c_{n,\xi}g_n}\bigr)
 \le \|b_n(\xi)\|^2\,
 \mathcal E_n\bigl(V_{e^{-t_n\|b_n\|^2/2}b_ng_n}\bigr)
 \longrightarrow0.
\]
The projected fields are scalar tangent Helmholtz fields.
Their harmonic traces
\[
 h_{n,\xi}=e^{-t_n\|b_n\|^2/2}c_{n,\xi}g_n
       \longrightarrow c_\xi g\quad\hbox{in }W^{-1,r}
\]
have uniformly bounded $L^r$ representatives, and their Poisson
extensions converge strongly on compact subdisks by
Lemma~\ref{lem:allp2-smoothing}.

Fix $\xi$. Let $a_{k,n}\in\mathbb R$ and $Y_n$ be the coefficients
and remainder in the Hodge decomposition
\[
 V_{h_{n,\xi}}\circ\Phi_n^{-1}
 =\sum_{k=1}^m a_{k,n}e_{k,n}+Y_n,
 \qquad Y_n\perp e_{k,n}\quad(1\le k\le m).
\]
By \eqref{eq:sc-residual-gap},
$\|Y_n\|_{L^2(\Omega_{\rho_n})}\to0$. The local boundedness of
$Ah_{n,\xi}$ and \eqref{eq:sc-mass-coercivity} bound the coefficients.
After passing to a subsequence, let $a_k=\lim_n a_{k,n}$.
For each disk $B\Subset\mathbb D$, the remainder estimate in the proof
of Lemma~\ref{lem:sc-point-evaluation} gives a constant $C_B$ independent
of $n$ such that
\[
 \|(I-T_n)(Y_n\circ\Phi_n)\|_{L^2(B)}
 \le C_B\|Y_n\|_{L^2(\Omega_{\rho_n})}\longrightarrow0.
\]
Consequently,
\begin{align*}
 Ah_{n,\xi}
 &=(I-T_n)V_{h_{n,\xi}}\\
 &=\sum_{k=1}^m a_{k,n}L_{k,n}
       +(I-T_n)(Y_n\circ\Phi_n)\\
 &\longrightarrow \frac1{\sqrt\mu}\sum_{k=1}^m a_kH_k
 \quad\hbox{in }L^2(B).
\end{align*}
Define $u_\xi=\mu^{-1/2}\sum_{k=1}^m a_ku_k\in E_\mu$.
The last limit is $Ag_{u_\xi}$, while the already established local
convergence of the Poisson extensions gives $A(c_\xi g)$.
Injectivity of the Poisson extension therefore yields
\begin{equation}\label{eq:allp2-multiplier-traces}
 g_{u_\xi}=c_\xi g.
\end{equation}
This constructs the limiting field without assuming a finite-energy
extension for an arbitrary $W^{-1,r}$ trace.
Separate subsequences for different $\xi$ cause no ambiguity because
the map $u\mapsto g_u$ is injective by
Lemma~\ref{lem:allp2-support-not-point}.

The functions $c_\xi$ separate points of $\mathbb S^1$. If their
values agree at $s,t$, then $b(s)-b(t)$ is orthogonal to the closed
span of $b(\xi)$ for the chosen dense set. By continuity, the
difference also belongs to that span and must therefore be zero.
Injectivity of $b$ gives $s=t$.

\section{Interior critical points and global coordinates}
\subsection{Vanishing trace on a boundary arc}
\label{sec:finite-area-missing-arc}

We first prove a radial estimate under the finite-area condition
$F\in A^2$. It will show that a pulled-back eigenfunction is
constant on any boundary arc where its harmonic gradient trace vanishes.

\begin{lemma}[A radial product estimate]
\label{lem:finite-area-radial-product}
Let $F\in A^2(\mathbb D)$ and $Z\in H^1_0(\mathbb D)$.
For $\rho\in(1/2,1)$, put
$A_\rho=\{z:\rho<|z|<1\}$.
Then, as $\rho\uparrow1$,
\begin{equation}\label{eq:finite-area-product-tail}
 \|F(\rho\,\cdot) Z(\rho\,\cdot)\|_{L^1(\mathbb S^1)}
 \le
 \left(\frac{2\log(1/\rho)}{1-\rho^2}\right)^{1/2}
 \|F\|_{L^2(A_\rho)}\|\nabla Z\|_{L^2(A_\rho)}
 \longrightarrow0.
\end{equation}
The radial slices of $Z$ are understood in their continuous
$L^2(\mathbb S^1)$ representative on every closed outer annulus.
\end{lemma}

\begin{proof}
Since $F$ is analytic, $t\mapsto\|F(t\,\cdot)\|_{L^2(\mathbb S^1)}^2$
is nondecreasing. Hence
\[
 \|F\|_{L^2(A_\rho)}^2
 =\int_\rho^1\|F(t\,\cdot)\|_{L^2(\mathbb S^1)}^2t\,dt
 \ge\frac{1-\rho^2}{2}\|F(\rho\,\cdot)\|_{L^2(\mathbb S^1)}^2.
\]
On the annulus, the map $\rho\mapsto Z(\rho\,\cdot)$ belongs to
$H^1((1/2,1);L^2(\mathbb S^1))$ and has zero trace at $\rho=1$.
The Hilbert-valued fundamental theorem of calculus and Cauchy--Schwarz
therefore give
\[
 \begin{aligned}
 \|Z(\rho\,\cdot)\|_2
 &\le\int_\rho^1\|\partial_t[Z(t\,\cdot)]\|_2\,dt\\
 &\le\left(\int_\rho^1\frac{dt}{t}\right)^{1/2}
      \left(\int_\rho^1t\|\partial_t[Z(t\,\cdot)]\|_2^2\,dt\right)^{1/2}\\
 &\le\sqrt{\log(1/\rho)}\,\|\nabla Z\|_{L^2(A_\rho)}.
 \end{aligned}
\]
Multiply these estimates and apply Cauchy--Schwarz on $\mathbb S^1$.
The prefactor in \eqref{eq:finite-area-product-tail} remains bounded
as $\rho\uparrow1$ (in fact it tends to one). Both norms over $A_\rho$ tend
to zero by absolute continuity of their integrals.
\end{proof}

\begin{lemma}[Vanishing trace on an arc]
\label{lem:finite-area-missing-arc}
Let $F$ be zero-free and belong to $A^2(\mathbb D)$, let
$v\in H^1(\mathbb D;\mathbb R)$, and suppose
\[
 W=\frac{2v_{\bar z}}{\overline F}=H+Z,
 \qquad Z\in H^1_0(\mathbb D),\qquad H=Ag,
 \qquad g\in\mathcal D'(\mathbb S^1;\mathbb C),
\]
where $A$ is disk Poisson extension. If $g$ vanishes on an open
connected circle arc $I$, then the Sobolev trace of $v$ is constant
on $I$. If $v$ is continuous on the closed disk, the boundary value
is pointwise constant on $I$.
\end{lemma}

\begin{proof}
Fix a compact subarc $I'\Subset I$. Choose a smooth cutoff supported
in $I$ and identically one near $I'$. Its product with $g$ is zero.
Inserting the complementary cutoff into the distributional Poisson
representation shows that $H$ is smooth up to $I'$ and vanishes
there. Indeed, the source support stays separated from $I'$, so all
source derivatives of the Poisson kernel have the factor $1-|z|^2$
and a denominator bounded away from zero. Thus
\[
 \sup_{e^{i\theta}\in I'}|H(\rho e^{i\theta})|\le C_{I'}(1-\rho)
\]
for $\rho$ sufficiently close to one.

For almost every $\rho$, the chain rule gives
\[
 |\partial_\theta v(\rho e^{i\theta})|
 \le \rho|\nabla v(\rho e^{i\theta})|
 =\rho|F(\rho e^{i\theta})|\,|H(\rho e^{i\theta})+Z(\rho e^{i\theta})|.
\]
The $F(\rho\,\cdot)Z(\rho\,\cdot)$ term tends to zero in $L^1(\mathbb S^1)$ by
Lemma~\ref{lem:finite-area-radial-product}. For the other term,
\[
 \begin{aligned}
 \|F(\rho\,\cdot)H(\rho\,\cdot)\|_{L^1(I')}
 &\le C_{I'}(1-\rho)\|F(\rho\,\cdot)\|_{L^2(\mathbb S^1)}\\
 &\le C_{I'}\sqrt{1-\rho}\,\|F\|_{L^2(A_\rho)}
 \longrightarrow0.
 \end{aligned}
\]
Consequently
\begin{equation}\label{eq:finite-area-tangential-variation}
 \|\partial_\theta v(\rho\,\cdot)\|_{L^1(I')}\longrightarrow0
\end{equation}
along the full-measure set of admissible radii.

The annular $H^1$ regularity of $v$ gives
$v(\rho\,\cdot)\to\operatorname{Tr}v$ in $L^2(\mathbb S^1)$ as $\rho\uparrow1$.
For every $\psi\in C_c^\infty(I')$,
\[
 \int_{I'}(\operatorname{Tr}v)\psi'
 =\lim_{\rho\uparrow1}\int_{I'}v(\rho\,\cdot)\psi'
 =-\lim_{\rho\uparrow1}\int_{I'}(\partial_\theta v(\rho\,\cdot))\psi=0.
\]
Thus the distributional tangential derivative of the trace is zero
on $I$. Since $I$ is connected, its trace is constant there.
Continuity, when available, gives the pointwise assertion.
\end{proof}

\subsection{The support of the gradient trace}
\label{sec:allp2-support}

\begin{lemma}\label{lem:allp2-support-not-point}
Under the conformal hypothesis of Section~\ref{sec:sa-spectral},
for every nonzero real $u\in E_\mu$, the harmonic gradient trace
$g_u$ defined in Section~\ref{sec:scalar-multiplier-limits}
is nonzero and cannot be supported at a single point.
The linear map $u\mapsto g_u$ on $E_\mu$ is injective.
\end{lemma}
\begin{proof}
With $v=u\circ\Phi$ and $W=2v_{\bar z}/\overline F$, write
$H=(I-T)W=Ag_u$ for the harmonic remainder.
If $g_u=0$, then $H=0$, hence $W=TW$ in $L^2(q)$.
The strict Dirichlet form bound gives $\|T\|<1$, so $W=0$ and
$u=0$, a contradiction. This also proves injectivity by taking
differences.

The identity $W_z=-\mu Fv/2\in L^2$ and $TW\in H^1_0$ imply
$H_z\in L^2$. Write $\widehat{g_u}(n)$ for the coefficient of
$e^{in\theta}$ in the distributional Fourier series of $g_u$.
Harmonic Fourier orthogonality gives
\[
 \sum_{n\ge1}n|\widehat{g_u}(n)|^2<\infty.
\]
Thus the positive Fourier coefficients tend to zero. A distribution
supported at one point is a finite sum of derivatives of a Dirac
mass; its positive Fourier coefficients are an exponential factor
times a polynomial. They cannot tend to zero unless the distribution
vanishes. The decay follows from the identity for $W_z$;
membership in $W^{-1,r}$ alone would not exclude a Dirac mass.
\end{proof}

\begin{proof}[Proof of Theorem~\ref{thm:allp2-main}]
Suppose that a nonzero real $u\in E_\mu$ has an interior critical
point, and put $g=g_u$. Section~\ref{sec:scalar-multiplier-limits}
provides functions $c_\xi$ separating the points of $\mathbb S^1$ and
eigenfunctions $u_\xi$ satisfying \eqref{eq:allp2-multiplier-traces}.
By linearity, every $c\in\operatorname{span}_{\mathbb R}\{1,c_\xi\}$
gives the trace $cg$ of an eigenfunction in $E_\mu$. That
eigenfunction is unique by Lemma~\ref{lem:allp2-support-not-point}.

If $\operatorname{supp}g=\mathbb S^1$, the local inversion argument in
Section~\ref{sec:allp2-smoothing} shows that $cg=0$ forces $c=0$.
Thus this correspondence is injective, and
Lemma~\ref{lem:sa-multiplicity} gives
$\dim\operatorname{span}_{\mathbb R}\{1,c_\xi\}\le2$.
At least one $c_\xi$ is nonconstant, since the functions separate
points. All the others are affine functions of it. It would
therefore define a continuous injection of the circle into
$\mathbb R$, which is impossible.

Otherwise $g$ vanishes on a nonempty open arc $I$. There is a $\xi$
for which $u,u_\xi$ are linearly independent. Otherwise write
$u_\xi=\lambda_\xi u$ for each $\xi$. Then
$(c_\xi-\lambda_\xi)g=0$ for every $\xi$, and local inversion
places the support of $g$ in a common level set of all the $c_\xi$.
Such a set contains at most one point, contradicting
Lemma~\ref{lem:allp2-support-not-point}.

Both traces $g$ and $c_\xi g$ vanish on $I$. Locality of the $C^1$
multiplication from
Section~\ref{sec:allp2-smoothing} follows by approximating compactly
supported $W^{1,2p/(p-1)}$ tests on a smaller arc by smooth tests.
Lemma~\ref{lem:finite-area-missing-arc} applies with $F\in A^2$
and $Z=TW\in H^1_0$; thus the two
pulled-back eigenfunctions are constant on $I$. Choose a nonzero
real linear combination $w$ of these pullbacks whose Dirichlet
trace is zero on $I$. It satisfies $-\Delta w=qw$ in $\mathbb D$,
and its weak Neumann data also vanish.

Let $\widetilde w$ and $\widetilde q$ be the zero extensions of
$w$ and $q$ across a smaller circle arc. The zero Dirichlet trace
gives $\widetilde w\in H^1_{\rm loc}$, and the weak Neumann
equation gives
\[
 -\Delta\widetilde w=\widetilde q\,\widetilde w
\]
with no distribution supported on the boundary arc. Since $w$ is
bounded and $q\in L^p$, $p>1$, interior elliptic
regularity gives $\widetilde w\in W^{2,p}_{\rm loc}$. Moreover,
\[
 \sup_x\int_{|x-y|<R}|\log|x-y||\,|\widetilde q(y)|\,dy
 \le\|\widetilde q\|_p
 \left(\int_{|z|<R}|\log|z||^{p'}\,dA(z)\right)^{1/p'}
 \longrightarrow0\qquad(R\downarrow0).
\]
Sawyer's planar unique-continuation theorem \cite{Sawyer1984}
therefore applies to this local Kato potential. The extension is
zero on the outside open set, so $w$ vanishes near the arc and then
on the connected disk. This contradicts the linear independence of
$u,u_\xi$, proving
that $u$ has no interior critical point.

Lemma~\ref{lem:sa-conformal} gives
$v=u\circ\Phi\in C(\overline{\mathbb D})$.
Since $\Phi$ is conformal, the chain rule shows that
$\nabla v$ is nonzero in $\mathbb D$.
Both extrema of $v$ on $\overline{\mathbb D}$ are attained,
and an interior extremum would have zero gradient.
Thus both extrema are attained only on $\mathbb S^1$,
which proves \eqref{eq:intro-prime-end-hotspots}.
Compactness and the equivalence of the domain hypotheses follow from
Lemma~\ref{lem:sa-conformal} and
Proposition~\ref{prop:qhb-conformal-regularity}.
\end{proof}
\Needspace{10\baselineskip}
\subsection{Global coordinates and the stream function}
\label{subsec:global-eigen-coordinates}

\begin{corollary}[Global eigenfunction coordinates]
\label{cor:local-eigen-coordinates}
Let $\Omega$ satisfy the hypotheses of Theorem~\ref{thm:allp2-main}.
If $\dim_{\mathbb R}E_\mu=2$ and $u_1,u_2$ form a real basis, then
\[
 \mathcal U:\Omega\longrightarrow\mathbb R^2,
 \qquad \mathcal U(x)=(u_1(x),u_2(x)),
\]
is a real-analytic diffeomorphism onto a bounded domain $\Sigma$
star-shaped about $0$. Its Jacobian has constant nonzero sign,
and there is exactly one $x_*\in\Omega$ with
$u_1(x_*)=u_2(x_*)=0$.
\end{corollary}
\begin{proof}
If the two gradients were dependent at a point, a nonzero real
linear combination would have a critical point there, contrary to
Theorem~\ref{thm:allp2-main}. Thus $\det D\mathcal U$ is nowhere
zero and has constant sign.

Fix a conformal bijection $\Phi:\mathbb D\to\Omega$ with
$F=\Phi'\in L^{2p}(\mathbb D)$, $p>1$, and put
$v_k=u_k\circ\Phi$ for $k=1,2$ and
$\widehat{\mathcal U}=(v_1,v_2)$.
After interchanging the basis elements if necessary, we may assume
$j=\det D\widehat{\mathcal U}>0$ in $\mathbb D$.
The components are real analytic in $\mathbb D$.
Lemma~\ref{lem:sa-conformal} gives
$v_k\in W^{2,p}(\mathbb D)\cap C(\overline{\mathbb D})$ and
$-\Delta v_k=\mu|F|^2v_k$ with homogeneous Neumann data.
In particular, $\nabla v_k\in L^{2p}(\mathbb D)$ and
$j\in L^p(\mathbb D)$.

The common equation makes the system
\[
 \psi_x=v_2(v_1)_y-v_1(v_2)_y,\qquad
 \psi_y=v_1(v_2)_x-v_2(v_1)_x
\]
compatible. Its right side belongs to $W^{1,p}$, so the Poincar\'e
lemma on the disk gives $\psi\in W^{2,p}(\mathbb D)$, unique up to
a constant. Its tangential boundary derivative is
$v_1\partial_\nu v_2-v_2\partial_\nu v_1=0$, where $\nu$ is the
outward unit normal to $\mathbb S^1$. These product traces follow by
approximation in $W^{2,p}$, using uniform convergence of the
undifferentiated factors. Thus the boundary value of $\psi$ is
constant; normalize it to zero. Direct differentiation gives
\begin{equation}\label{eq:eigen-stream}
 -\Delta\psi=2j,\qquad
 D\widehat{\mathcal U}\,\nabla\psi=-j\widehat{\mathcal U},\qquad
 \psi|_{\mathbb S^1}=0.
\end{equation}
The maximum principle gives $\psi>0$ in $\mathbb D$.
Since $\psi$ is continuous on the closed disk, each closed
superlevel set $\{\psi\ge c\}$, $c>0$, is compact in $\mathbb D$.

Choose a maximizer $z_*\in\mathbb D$ of $\psi$.
Equation~\eqref{eq:eigen-stream} gives
$\widehat{\mathcal U}(z_*)=0$.
We apply the lifted-ray criterion of Gorni and Zampieri
\cite[Theorem~1.1 and Proposition~1.2]{GorniZampieri1994}.
The function $-\log\psi$ is smooth, proper on $\mathbb D$, and
bounded below: its closed sublevel sets are compact positive
superlevel sets of $\psi$, and $\psi$ is bounded above. Moreover,
\[
 \nabla(-\log\psi)\cdot[D\widehat{\mathcal U}]^{-1}
           \widehat{\mathcal U}
 =\frac{|\nabla\psi|^2}{j\psi}\ge0.
\]
Since $\mathbb D$ is open and connected and
$\widehat{\mathcal U}$ is a local diffeomorphism, the criterion
gives injectivity and an image star-shaped about
$\widehat{\mathcal U}(z_*)=0$.
The local analytic inverses give a global analytic inverse, and
continuity on the closed disk gives boundedness of the image.
Its unique zero corresponds to $x_*=\Phi(z_*)$ in $\Omega$.
\end{proof}

Once $j$ has fixed nonzero sign, the argument on the disk uses only
interior smoothness, the common scalar equation, homogeneous Neumann
data, and $W^{2,p}$ regularity with $p>1$; nodal connectedness is not a
separate input.

For $z_0\in\mathbb D$, let $z(t)$ solve
$\dot z=\nabla\psi(z)/j(z)$ with $z(0)=z_0$.
The stream identities give
\begin{equation}\label{eq:eigen-radial-flow}
 \widehat{\mathcal U}(z(t))=e^{-t}\widehat{\mathcal U}(z_0),\qquad
 \frac{d}{dt}\psi(z(t))=\frac{|\nabla\psi(z(t))|^2}{j(z(t))}\ge0.
\end{equation}
The trajectory stays in
$\{\psi\ge\psi(z_0)\}\Subset\mathbb D$, so it exists for all
$t\ge0$.
Equation~\eqref{eq:eigen-stream} also shows that $z_*$ is the only
critical point of $\psi$. Differentiating its second identity at
$z_*$, where $\nabla\psi=0$ and $D\widehat{\mathcal U}$ is
invertible, shows that the Hessian of $\psi$ is $-j(z_*)$ times
the identity matrix.
Away from $z_*$, the gradient trajectories map to radial segments.
Along the preimage of each ray from $0$, $\psi$ decreases strictly as the radius
increases and tends to zero at the image boundary: the inverse
points leave every compact subset of the disk, and $\psi$ tends
to zero at the circle. Thus, for $0<c<\psi(z_*)$, the level set
$\{\psi=c\}$ is a simple closed curve meeting each such preimage
once and orthogonally. The simple-curve assertion follows by the
implicit function theorem applied to its radius as a function of
the target angle.

\Needspace{10\baselineskip}
\subsection{The boundary angle and the coordinate image}
\label{subsec:eigen-coordinate-boundary}

\begin{proposition}[Conformal boundary coordinates]
\label{prop:eigen-coordinate-boundary}
Under the hypotheses of Corollary~\ref{cor:local-eigen-coordinates},
fix a conformal bijection $\Phi:\mathbb D\to\Omega$ with
$\Phi'\in L^{2p}(\mathbb D)$, $p>1$.
Put $\widehat{\mathcal U}=\mathcal U\circ\Phi$ and order the basis
so that $j=\det D\widehat{\mathcal U}>0$.
The map $\widehat{\mathcal U}$ extends to a homeomorphism from
$\overline{\mathbb D}$ onto $\overline\Sigma$.
Identify the target $\mathbb R^2$ with $\mathbb C$.
The boundary trace $\gamma(t)=\widehat{\mathcal U}(e^{it})$,
$t\in\mathbb R$, is absolutely continuous and nonzero, and has
an argument lift $\alpha$ with
\[
 \gamma(t)=|\gamma(t)|e^{i\alpha(t)},\qquad
 \alpha(t+2\pi)=\alpha(t)+2\pi.
\]
The induced angular map is a $C^{1,\beta}$ diffeomorphism of the
circle, where one may take
\begin{equation}\label{eq:eigen-angle-holder}
 \beta=2-\frac2p\quad(1<p<2),\qquad
 0<\beta<1\quad(p\ge2).
\end{equation}
For every $t\in\mathbb R$,
\begin{equation}\label{eq:eigen-angle-poisson}
 |\gamma(t)|^2\alpha'(t)
 =\frac1\pi\int_{\mathbb D}
    \frac{1-|z|^2}{|e^{it}-z|^2}\,j(z)\,dA(z).
\end{equation}
This continuous product agrees with
$\det(\gamma(t),\gamma'(t))$ almost everywhere. Moreover,
\begin{equation}\label{eq:eigen-angle-lower}
 \alpha'(t)\ge
 \frac1{3\pi\|\gamma\|_\infty^2}
       \int_{|z|\le1/2}j(z)\,dA(z)>0
 \qquad(t\in\mathbb R).
\end{equation}
The inverse boundary map $\partial\Sigma\to\mathbb S^1$ is
Lipschitz for Euclidean distance. There is a continuous positive
$2\pi$-periodic function $R:\mathbb R\to(0,\infty)$ such that
\begin{equation}\label{eq:eigen-image-radial}
 \Sigma=\{\rho(\cos\theta,\sin\theta):
       \theta\in\mathbb R,\ 0\le\rho<R(\theta)\}.
\end{equation}
For a Jordan domain $\Omega$, $\mathcal U$ extends to a
homeomorphism of Euclidean closures; in general the source closure
is the conformal, or prime-end, compactification.
\end{proposition}
\begin{proof}
Lemma~\ref{lem:sa-conformal} gives continuity of
$\widehat{\mathcal U}$ on the closed disk. No boundary value lies
in $\Sigma$: otherwise interior points tending to that boundary
point would have images converging in $\Sigma$, contradicting
continuity of the inverse from Corollary~\ref{cor:local-eigen-coordinates}.
Compactness then gives
$\widehat{\mathcal U}(\mathbb S^1)=\partial\Sigma$, which avoids $0$.
The unique interior zero and $j>0$ give degree one at $0$, hence
the stated lift $\alpha$. The trace theorem gives
$\gamma\in W^{2-1/p,p}(\mathbb S^1)$, so $\gamma$ and its local
argument lifts are absolutely continuous.

Let $\psi$ be the normalized stream function constructed in the
proof of Corollary~\ref{cor:local-eigen-coordinates}.
Its defining derivatives give, almost everywhere on $\mathbb S^1$,
\[
 -\partial_\nu\psi
 =\det(\gamma,\gamma')=|\gamma|^2\alpha'.
\]
For smooth real $h$ on $\mathbb S^1$, its Poisson extension $Ah$
satisfies, by Green's formula and $\psi|_{\mathbb S^1}=0$,
\[
 \int_0^{2\pi}h(e^{it})|\gamma(t)|^2\alpha'(t)\,dt
 =2\int_{\mathbb D}(Ah)(z)j(z)\,dA(z).
\]
Inserting the disk Poisson formula proves
\eqref{eq:eigen-angle-poisson} almost everywhere. Fubini's theorem
applies since $j\in L^1$ and
$\int_{\mathbb S^1}P_s(z)\,ds=1$.

To obtain boundary regularity, choose a closed collar of
$\mathbb S^1$ on which $|\widehat{\mathcal U}|$ is bounded away
from zero. Sobolev embedding gives
$|\widehat{\mathcal U}|^{-2}\in C^{0,\beta}$ there for every
exponent in \eqref{eq:eigen-angle-holder}.
For a local branch $\vartheta$ of $\arg\widehat{\mathcal U}$,
the stream identities give
\[
 \nabla\vartheta
 =|\widehat{\mathcal U}|^{-2}(\psi_y,-\psi_x),\qquad
 \operatorname{div}\bigl(|\widehat{\mathcal U}|^{-2}
                         \nabla\psi\bigr)=0.
\]
This is a uniformly elliptic divergence-form equation for
$\psi\in H^1$ with zero trace on the smooth outer boundary.
The boundary Schauder estimate for weak solutions
\cite[Corollary~8.36]{GilbargTrudinger2001} gives
$\psi\in C^{1,\beta}$ up to $\mathbb S^1$ on a smaller collar.
The first identity gives the same regularity for $\vartheta$.
In particular, $\alpha'$ is $C^{0,\beta}$. Since $\psi$ is
constant on $\mathbb S^1$, the same identity also gives
$\partial_\nu\vartheta=0$ there.

The integral in \eqref{eq:eigen-angle-poisson} is continuous in
$t$. Indeed, $\nabla|\widehat{\mathcal U}|\in L^{2p}$ in the
collar, while $\nabla\vartheta$ is bounded, so
\[
 j=|\widehat{\mathcal U}|
      \det\bigl(\nabla|\widehat{\mathcal U}|,
                \nabla\vartheta\bigr)\in L^{2p}
\]
on the collar. Interior smoothness then gives
$j\in L^{2p}(\mathbb D)$.
The kernel in \eqref{eq:eigen-angle-poisson} is bounded by
$2/|e^{it}-z|$, hence belongs to
$L^{2p/(2p-1)}(\mathbb D)$. Its dependence on $t$ is continuous
in that space by continuity of rotations. H\"older's inequality
therefore makes the integral continuous, proving the identity at
every $t$. On $|z|\le1/2$ the kernel is at least $1/3$, which
gives \eqref{eq:eigen-angle-lower}.

The positive derivative and degree one make the angular map a
$C^{1,\beta}$ circle diffeomorphism. Thus $\gamma$ is an injective
polar parametrization and defines the continuous radius $R$.
Its inverse is the composition of radial projection with the
inverse angular map. Radial projection is Lipschitz where the
radius is bounded away from zero, and the inverse angular map is
$C^1$ on the compact circle, hence Lipschitz also for chord distance.
The interior and boundary images are disjoint, and both restrictions
are bijective; compactness gives the homeomorphism of closures.
For a Jordan source, Carath\'eodory's theorem identifies the
conformal closed disk with its Euclidean closure.
\end{proof}
\begin{corollary}[Boundary zeros and nodal endpoints]
\label{cor:eigen-boundary-zeros}
Under the hypotheses of Proposition~\ref{prop:eigen-coordinate-boundary},
every nonzero $w\in E_\mu$ has exactly two prime-end boundary zeros.
Both are simple in the conformal parameter: the derivative of
$t\mapsto(w\circ\Phi)(e^{it})$ exists and is nonzero at each zero.
The closed nodal set of $w\circ\Phi$ is a single arc in
$\overline{\mathbb D}$. For every $\beta$ allowed by
\eqref{eq:eigen-angle-holder}, this arc is $C^{1,\beta}$ up to
its two endpoints and meets $\mathbb S^1$ orthogonally.
\end{corollary}
\begin{proof}
After multiplying $w$ by a nonzero constant, its boundary trace has
the form $|\gamma(t)|\sin(\alpha(t)-\theta_0)$ for some $\theta_0\in\mathbb R$.
The angular diffeomorphism gives exactly two zeros on the circle.
At either zero $t_0$, continuity of $|\gamma|$ and differentiability
of $\alpha$ give
\[
 \lim_{t\to t_0}
 \frac{|\gamma(t)|\sin(\alpha(t)-\theta_0)}{t-t_0}
 =|\gamma(t_0)|\cos(\alpha(t_0)-\theta_0)\alpha'(t_0)\ne0.
\]
The cosine equals $1$ or $-1$, and
\eqref{eq:eigen-angle-lower} makes the last factor positive.

The homeomorphism in Proposition~\ref{prop:eigen-coordinate-boundary}
identifies the closed nodal set with a segment through $0$ in
$\overline\Sigma$. Near either endpoint, it is a level set of a
local branch $\vartheta$ of $\arg\widehat{\mathcal U}$.
The proof of that proposition gives $\vartheta\in C^{1,\beta}$
up to the circle and
\[
 \nabla\psi=|\widehat{\mathcal U}|^2
        (-\partial_y\vartheta,\partial_x\vartheta).
\]
Since $\psi$ is constant on $\mathbb S^1$, this identity gives
$\partial_\nu\vartheta=0$ there, where $\nu$ is the outward unit
normal. The tangential derivative is $\alpha'>0$.
The implicit function theorem therefore gives a $C^{1,\beta}$
nodal arc whose tangent at the endpoint is normal to $\mathbb S^1$.
\end{proof}

\begin{corollary}[Regularity of the radial function]
\label{cor:eigen-radial-regularity}
Under the hypotheses and notation of
Proposition~\ref{prop:eigen-coordinate-boundary}, let $R$ be the
positive periodic radial function of $\partial\Sigma$.
Then
\[
 \begin{array}{ll}
 R\in W^{1,p/(2-p)}(\mathbb S^1), & 1<p<2,\\[2pt]
 R\in W^{1,r}(\mathbb S^1)\text{ for every }1\le r<\infty,
     & p=2,\\[2pt]
 R\in C^{1,1-2/p}(\mathbb S^1), & p>2.
 \end{array}
\]
In all cases $\partial\Sigma$ is rectifiable.
\end{corollary}
\begin{proof}
For $1<p<2$, the trace theorem and the one-dimensional Sobolev
embedding give
\[
 \gamma'\in W^{1-1/p,p}(\mathbb S^1)
       \subset L^{p/(2-p)}(\mathbb S^1).
\]
When $p=2$, the same trace space embeds into every finite
$L^r(\mathbb S^1)$. The function $|\gamma|$ is absolutely continuous,
with $|(|\gamma|)'|\le|\gamma'|$ almost everywhere. Since
$\alpha^{-1}$ is increasing and Lipschitz,
$R=|\gamma|\circ\alpha^{-1}$ is absolutely continuous.
Take $r=p/(2-p)$ when $1<p<2$, and any $1\le r<\infty$
when $p=2$. Periodicity and change of variables give
\[
 \begin{aligned}
 \int_0^{2\pi}|R'(\theta)|^r\,d\theta
 &=\int_0^{2\pi}|(|\gamma|)'(t)|^r
                   \alpha'(t)^{1-r}\,dt\\
 &\le\bigl(\min_{\mathbb R}\alpha'\bigr)^{1-r}
          \int_0^{2\pi}|\gamma'(t)|^r\,dt<\infty.
 \end{aligned}
\]
The minimum is positive by \eqref{eq:eigen-angle-lower}.

If $p>2$, Morrey's inequality gives
$\gamma\in C^{1,1-2/p}$. By \eqref{eq:eigen-angle-holder} and
\eqref{eq:eigen-angle-lower}, the angular map and its inverse have
at least this regularity. Thus $R=|\gamma|\circ\alpha^{-1}$
belongs to $C^{1,1-2/p}$.
In every case $R\in W^{1,1}(\mathbb S^1)$, so the parametrization
$\theta\mapsto R(\theta)e^{i\theta}$ is absolutely continuous
and has finite length.
\end{proof}

\begin{remark}[Smooth source domains]
In the setting of Proposition~\ref{prop:eigen-coordinate-boundary},
suppose in addition that $\partial\Omega$ is $C^\infty$.
The same conclusions can then be expressed in physical arclength.
Choose the basis orientation so
that $\det D\mathcal U>0$, let $s$ increase counterclockwise along
$\partial\Omega$, and let $\nu$ be the outward unit normal.
The stream construction and conformality give
\[
 \begin{aligned}
 -\Delta(\psi\circ\Phi^{-1})&=2\det D\mathcal U,\\
 \nabla(\psi\circ\Phi^{-1})
    &=u_1\nabla^\perp u_2-u_2\nabla^\perp u_1.
 \end{aligned}
\]
Here $\nabla^\perp h=(-h_y,h_x)$, as in Lemma~\ref{lem:sa-hodge}.
The function $\psi\circ\Phi^{-1}$ is positive in $\Omega$ and
vanishes on $\partial\Omega$. Hopf's boundary lemma therefore gives
\[
 \det(\mathcal U,\partial_s\mathcal U)
    =-\partial_\nu(\psi\circ\Phi^{-1})>0
 \qquad\text{on }\partial\Omega.
\]
The boundary eigenfunctions are smooth, so the angular map in
arclength is a smooth diffeomorphism and $R$ is smooth.
The two zeros in Corollary~\ref{cor:eigen-boundary-zeros} are simple
in physical arclength. At either zero, the nonzero tangential
derivative and the Neumann condition show that the nodal arc meets
$\partial\Omega$ orthogonally.
\end{remark}
\subsection{Convexity and star centers}
\label{subsec:eigen-coordinate-convexity}

A point $\xi$ is a \emph{star center} of $\Sigma$ if $\xi\in\Sigma$
and $[\xi,z]\subset\Sigma$ for every $z\in\Sigma$.
If $\Phi'\in A^{2p}$ for some $p>2$, the star centers contain a
neighborhood of $0$. Indeed, Corollary~\ref{cor:eigen-radial-regularity}
and the angular regularity give $\gamma\in C^1$, while
\eqref{eq:eigen-angle-poisson} and \eqref{eq:eigen-angle-lower} give
$\min_t\det(\gamma(t),\gamma'(t))>0$. For all sufficiently small
$|\xi|$, we have $\xi\in\Sigma$ and
\[
 \det(\gamma(t)-\xi,\gamma'(t))
 \ge \min_t\det(\gamma(t),\gamma'(t))
       -|\xi|\|\gamma'\|_\infty>0.
\]
The argument of $\gamma(t)-\xi$ therefore increases strictly and
has degree one.
Each ray from $\xi$ meets the boundary exactly once, so $\xi$ is a star
center. For a smooth source, the same proof uses physical arclength
and the strict angular derivative obtained above by Hopf's lemma.

\begin{proposition}\label{prop:eigen-image-convexity}
Under the hypotheses of Corollary~\ref{cor:local-eigen-coordinates},
$\Sigma$ is convex if and only if
\[
 \Omega\setminus\{w=c\}
 \quad\hbox{has at most two connected components}
 \qquad(0\ne w\in E_\mu,\ c\in\mathbb R).
\]
There are smooth strictly convex sources with a double first positive
Neumann eigenvalue for which $\Sigma$ is nonconvex for every real basis.
\end{proposition}
\begin{proof}
Each strict side of an affine line intersects a convex open set in a
connected set or the empty set. Pullback by $\mathcal U$ proves one
direction.

If $\Sigma$ is nonconvex, some line meets it in at least two disjoint
open intervals. Each interval is a properly embedded line in the
topological disk $\Sigma$. The first separates $\Sigma$ into two open
disks, by Jordan--Schoenflies after one-point compactification of
$\Sigma$. The second is proper in one of these disks and separates
it again. Each of the three resulting open regions contains points
off the affine line; removing the remaining intervals cannot reconnect
them. Pulling back gives at least three components of a shifted level
complement.

B\'erard and Helffer
\cite[Theorem~1.2 and Proposition~3.5]{BerardHelffer2021}
construct smooth strictly convex domains whose first positive Neumann
eigenspace is two-dimensional and for which some $w-c$, with
$w\in E_\mu$, has exactly three nodal domains. Combining their result
with Corollary~\ref{cor:local-eigen-coordinates} and the criterion just
proved gives nonconvex coordinate images. An invertible real linear
change of basis preserves convexity, proving the last assertion.
\end{proof}

\paragraph{The equilateral triangle.}
On the triangle with vertices $0$, $1$, and $e^{i\pi/3}$, the real
and imaginary parts of the function
\[
 Z(x,y)=\frac13\left(e^{-4\pi ix/3}
              +2e^{2\pi ix/3}\cos(2\pi y/\sqrt3)\right)
\]
form a basis of the first positive Neumann eigenspace.
The eigenvalue is $16\pi^2/9$ and is double
\cite[Section~3.1 and Appendix~A.2]{BerardHelffer2018}.
We take $\mathcal U=(\operatorname{Re}Z,\operatorname{Im}Z)$.
After rescaling the source and complex conjugation, this is the
diffeomorphism from the triangle onto the deltoid described by
Bakry and Zribi \cite[Section~3]{BakryZribi2016}. The example does
not depend on Theorem~\ref{thm:allp2-main}.
On the bottom side set $t=2\pi x/3$. Rotating gives the other sides
and the full boundary parametrization
\begin{equation}\label{eq:eigen-deltoid}
 d(t)=\frac{2e^{it}+e^{-2it}}3,
 \qquad 0\le t\le2\pi.
\end{equation}
It satisfies
\begin{equation}\label{eq:eigen-deltoid-angle}
 |d(t)|^2=\frac{5+4\cos3t}{9},\qquad
 \frac{d}{dt}\arg d(t)=\frac{2(1-\cos3t)}{5+4\cos3t}.
\end{equation}
The argument derivative is positive except at the three cusps, and
$d(t)=e^{it}(2+e^{-3it})/3$ has winding number one. These identities
also verify directly that the deltoid is star-shaped about $0$.
Put $\omega=e^{2\pi i/3}$. The midpoint $(1+\omega)/2$ has radius
$1/2$ and angle $\pi/3$, whereas $d(\pi/3)=e^{i\pi/3}/3$.
It lies outside the closed image. Since the closure of a convex set
is convex, the image is nonconvex.

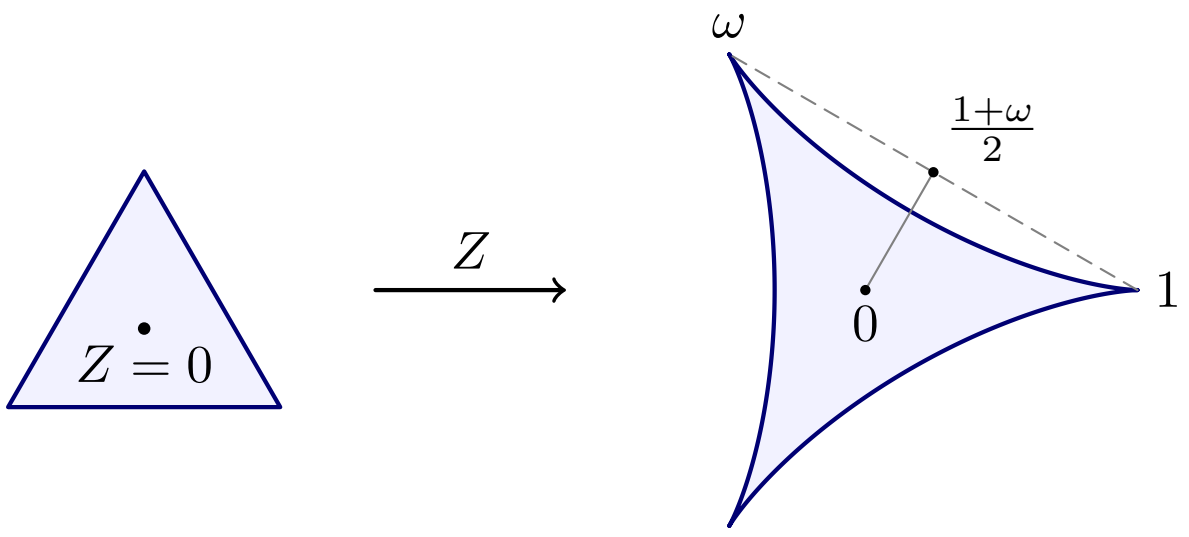
\begin{figure}[htbp]
\centering
\begin{tikzpicture}[x=1.8cm,y=1.8cm,line cap=round,line join=round,
                    every node/.style={font=\small}]
 \begin{scope}[shift={(-2,-.43)}]
  \filldraw[fill=blue!5,draw=blue!45!black,thick]
    (0,0)--(1,0)--(.5,{sqrt(3)/2})--cycle;
  \fill (.5,{sqrt(3)/6}) circle[radius=.023];
  \node[below,inner sep=3pt] at (.5,{sqrt(3)/6}) {$Z=0$};
 \end{scope}
 \draw[->,thick] (-.65,0)--(.05,0) node[midway,above] {$Z$};
 \begin{scope}[shift={(1.15,0)}]
  \draw[fill=blue!5,draw=blue!45!black,thick]
    plot[domain=0:360,samples=361,variable=\t]
      ({(2*cos(\t)+cos(2*\t))/3},
       {(2*sin(\t)-sin(2*\t))/3})--cycle;
  \draw[densely dashed,gray] (1,0)--(-.5,{sqrt(3)/2});
  \draw[gray] (0,0)--(.25,{sqrt(3)/4});
  \fill (.25,{sqrt(3)/4}) circle[radius=.019];
  \node[above right,inner sep=2pt] at (.25,{sqrt(3)/4})
       {$\frac{1+\omega}{2}$};
  \fill (0,0) circle[radius=.019];
  \node[below,inner sep=3pt] at (0,0) {$0$};
  \node[right,inner sep=3pt] at (1,0) {$1$};
  \node[above,inner sep=3pt] at (-.5,{sqrt(3)/2}) {$\omega$};
 \end{scope}
\end{tikzpicture}
\caption{The triangle-to-deltoid coordinates
\cite[Section~3]{BakryZribi2016}, with boundary given by
\eqref{eq:eigen-deltoid}. The midpoint of the dashed segment lies
outside the closed image.}
\label{fig:eigen-triangle-deltoid}
\end{figure}

At the cusp $d(0)=1$, using the argument branch with value zero there,
\[
 |d(t)|=1-t^2+O(t^4),\qquad
 (|d|)'(t)=-2t+O(t^3),\qquad
 \arg d(t)=t^3/3+O(t^5).
\]
Together with \eqref{eq:eigen-deltoid-angle}, these give
\[
 1-R(\theta)=3^{2/3}|\theta|^{2/3}+O(|\theta|^{4/3})
 \qquad(\theta\to0).
\]
For nonzero $\theta$ near zero, $|R'(\theta)|$ is bounded above and
below by positive multiples of $|\theta|^{-1/3}$.
The other cusps are rotations of this one, and $R$ is analytic
elsewhere. Hence $R\in W^{1,r}$ for $1\le r<3$, but
$R\notin W^{1,3}$.
If $\zeta\in\mathbb S^1$ maps to a triangle vertex, the
Schwarz--Christoffel formula
\cite[Sections~2.2 and~4.1]{DriscollTrefethen2002} bounds
$|\Phi'(z)|$ above and below by positive multiples of
$|z-\zeta|^{-2/3}$ near $\zeta$. Among exponents $p>1$, we have
$|\Phi'|^2\in L^p$ exactly when $p<3/2$; the area integral diverges logarithmically at
$p=3/2$. Corollary~\ref{cor:eigen-radial-regularity} recovers every
$r<3$ by letting $p$ increase to $3/2$.

To examine the boundary angular regularity, choose $\Phi$ so that
the preimages of the three vertices are the cube roots of unity and $\Phi(1)$ is the
triangle vertex $0$. For a nonzero complex constant $C$,
\[
 \Phi'(z)=C(1-z^3)^{-2/3},
\]
with the analytic branch in $\mathbb D$.
Let $\tau$ be the source circle angle near $1$, and use the signed
parameter $t$ of $d$ near the cusp $1$. On each adjacent side, $|t|$
is $2\pi/3$ times the distance to the vertex. Integrating
$|\Phi'(e^{i\tau})|=|C||1-e^{3i\tau}|^{-2/3}$ therefore gives
the following expansions for a constant $k>0$:
\[
 t=k\operatorname{sgn}(\tau)|\tau|^{1/3}(1+O(\tau^2)),\qquad
 \frac{dt}{d\tau}=\frac{k}{3}|\tau|^{-2/3}(1+O(\tau^2)).
\]
The exact formula \eqref{eq:eigen-deltoid-angle} gives
$(\arg d)'=t^2+\frac54t^4+O(t^6)$, so
\[
 \alpha'(\tau)=\frac{k^3}{3}
       +\frac{5k^5}{12}|\tau|^{2/3}+O(|\tau|^{4/3}).
\]
Thus the conformal boundary angle fails to be $C^{1,\beta}$ for
every $\beta>2/3$, although its derivative at the cusp is positive.
The thresholds $r=3$ and $\beta=2/3$ describe this example; these
calculations do not establish optimality of the general estimates
for each fixed $p$.

\Needspace{7\baselineskip}
The origin is the deltoid's only star center. At a regular boundary
point a star center $\xi$ must satisfy
$\det(d(t)-\xi,d'(t))\ge0$: the segment toward $\xi$ lies on the inward
side of the tangent to first order. At the cusp $1$,
\[
 \det(d(t)-\xi,d'(t))=-2t\operatorname{Im}\xi+O(t^2).
\]
Taking both signs of $t$ forces $\operatorname{Im}\xi=0$.
The rotated conditions at $\omega$ and $\omega^2$ then force $\xi=0$.
\appendix

\section{Reciprocal positivity and Hilbert realization}
\label{sec:rkchecked}

The regularization and reciprocal-kernel construction are due to
\cite[Section 5]{OAI}. We use the global matrix flow of Pineda Montoya,
Gallego Gonz\'alez, Porras Barrera and R\'ua Mu\~noz
\cite[Lemmas 4.4--4.5]{Pineda2026} in place of successive rank-one
updates, and include its proof. Reciprocals throughout this appendix
are entrywise, and undecorated matrix norms are Euclidean operator
norms. The Green kernel $\mathcal G_0$, the Poisson kernel
$P_s$, and $N_0(s,t)=1/(\pi|s-t|^2)$ have the normalizations of
Section~\ref{sec:holder-kernel-audit}; boundary integrals use arclength.
The kernels $K_q,N_q$ are defined by the nonnegative chains
\eqref{eq:hk-K}--\eqref{eq:hk-R}.

\begin{theorem}[Compact-support reciprocal theorem]
\label{thm:rkchecked-bounded}
Suppose $q\ge0$ is bounded, has compact support in $\mathbb D$, and
\begin{equation}\label{eq:rkchecked-subcritical}
 \int_{\mathbb D}q|\psi|^2\,dA
 \le d\int_{\mathbb D}|\nabla\psi|^2\,dA
 \quad(\psi\in H^1_0(\mathbb D)),\qquad 0<d<1.
\end{equation}
For each $z_0\in\mathbb D$, $K_q(z_0,\cdot)$ is positive and bounded,
and $N_q(s,t)$ is finite and positive for $s\ne t$. The kernel
\begin{equation}\label{eq:rkchecked-D}
 D_q(s,t)=\frac{K_q(z_0,s)K_q(z_0,t)}{N_q(s,t)}\quad(s\ne t),
 \qquad D_q(s,s)=0,
\end{equation}
is conditionally negative semidefinite: for every finite family of
points $s_i\in\mathbb S^1$ and coefficients $a_i\in\mathbb R$ with
$\sum_i a_i=0$,
$\sum_{i,j}a_i a_jD_q(s_i,s_j)\le0$.
\end{theorem}

\subsection{The operator bound and the initial regularization}

Put $T=G(q\,\cdot)$ on $L^2(q)=L^2(\mathbb D,q\,dA)$.
Testing the weak Dirichlet equation for $Tf$ and applying
\eqref{eq:rkchecked-subcritical} gives
\[
 \|Tf\|_{L^2(q)}^2\le d\|\nabla Tf\|_2^2
 =d\langle f,Tf\rangle_{L^2(q)}.
\]
Thus $T$ is self-adjoint, positivity preserving, and $\|T\|\le d$.
Write $S=(I-T)^{-1}$; its norm is at most $(1-d)^{-1}$.

For a real symmetric matrix $J$, write $\mathsf P(J)$ if $J$ has
at most one positive eigenvalue. This property is preserved by
congruence, principal submatrices, subtraction of a positive
semidefinite matrix, and limits in a fixed dimension. We need bounded
continuous positive kernels $C_\varepsilon\uparrow\mathcal G_0$
as $\varepsilon\downarrow0$,
including its infinite diagonal, such that every finite reciprocal
matrix of $C_\varepsilon$ has $\mathsf P$.

For $x,y\in\mathbb D$ and $t\ge0$, set
\begin{gather*}
 h(x,y)=\frac{|x-y|^2}{(1-|x|^2)(1-|y|^2)},\qquad
 \mathcal L(t)=\frac1{\log(1+1/t)}\ (t>0),\quad \mathcal L(0)=0,
 \\
 f(t)=t+\tfrac12-\mathcal L(t).
\end{gather*}
We recall the regularization from \cite[Section 5]{OAI}.
For $0<\sigma<1$, let $\mathbb E_\sigma$ be expectation for the
probability density proportional to $v^{\sigma-1}/(1+v)$ on
$(0,\infty)$, and put $\eta=(1+v)^{-1}$. Integration by parts and
rescaling give
\[
 \mathbb E_\sigma(1-\eta)=\sigma,\qquad
 t^{1-\sigma}(t+1)^\sigma
 =\mathbb E_\sigma\frac{t(t+1)}{t+\eta}
 =\mathbb E_\sigma\left(t+1-\eta-
                 \frac{\eta(1-\eta)}{t+\eta}\right).
\]
Since $\mathcal L(t)=\int_0^1t^{1-\sigma}(t+1)^\sigma\,d\sigma$, it follows that
\begin{equation}\label{eq:rkchecked-mixture}
 f(t)=\int_0^1\mathbb E_\sigma
          \frac{\eta(1-\eta)}{t+\eta}\,d\sigma,
 \qquad 0\le f(t)\le f(0)=\tfrac12.
\end{equation}
The identity at $t=0$ follows by dominated convergence.
The kernel
\[
 k(x,y)=\frac1{1+h(x,y)}
 =\frac{(1-|x|^2)(1-|y|^2)}{|1-x\overline y|^2}
\]
is positive semidefinite: expand $(1-x\overline y)^{-1}$ as a
Gram series, multiply by its conjugate, and use positive diagonal
congruence. For $0<\eta<1$, the expansion in entrywise powers
$(h+\eta)^{-1}=\sum_{n\ge0}(1-\eta)^n k^{n+1}$ is therefore
positive semidefinite by the Schur product theorem. The series
converges also on the diagonal. Equation~\eqref{eq:rkchecked-mixture}
proves that $f(h)$ is positive semidefinite.

Congruence by $\operatorname{diag}(1-|x_i|^2)$ transforms
$h(x_i,x_j)+1/2$ into
$\tfrac12(1+|x_i|^2)(1+|x_j|^2)-2x_i\cdot x_j$, which has
$\mathsf P$. Consequently the finite reciprocals of
\begin{equation}\label{eq:rkchecked-Ceps}
 C_\varepsilon(x,y)
 =\frac1{4\pi[\mathcal L(h(x,y))+\varepsilon f(h(x,y))]},
 \qquad 0<\varepsilon\le1,
\end{equation}
have $\mathsf P$. Indeed, their reciprocals equal
$4\pi[h+1/2-(1-\varepsilon)f(h)]$.
Also $\mathcal L+\varepsilon f=\varepsilon(h+1/2)+(1-\varepsilon)\mathcal L
\ge\varepsilon/2$, so
\[
 0<C_\varepsilon\le\mathcal G_0,\qquad
 C_\varepsilon\uparrow\mathcal G_0,\qquad
 C_\varepsilon(x,x)=(2\pi\varepsilon)^{-1}.
\]
No positive semidefiniteness of $C_\varepsilon$ itself is asserted.

\subsection{The finite resolvent by a global flow}

\begin{lemma}[Pineda Montoya et al.\ {\cite[Lemmas 4.4--4.5]{Pineda2026}}]
\label{lem:rkchecked-update}
For a matrix with strictly positive entries, let $^{\circ-1}$ denote
entrywise reciprocal. Let $C$ be real symmetric with strictly positive
entries and $\mathsf P(C^{\circ-1})$. If $W$ is a diagonal matrix of
the same size with nonnegative entries
and $\|W^{1/2}CW^{1/2}\|<1$, then
\[
 M=C+CW^{1/2}(I-W^{1/2}CW^{1/2})^{-1}W^{1/2}C
\]
has strictly positive entries and $\mathsf P(M^{\circ-1})$.
\end{lemma}
\begin{proof}
Put $B=W^{1/2}CW^{1/2}$ and
\[
 M(t)=C+tCW^{1/2}(I-tB)^{-1}W^{1/2}C,\qquad 0\le t\le1.
\]
The norm bound makes the inverse well defined throughout this
interval. Its geometric series has nonnegative entries, so $M(t)$
has strictly positive entries. The identity
$M(t)W^{1/2}=CW^{1/2}(I-tB)^{-1}$ gives $M'=MWM$.
Writing $J=M^{\circ-1}$ and $W=\operatorname{diag}(w_k)$ yields
the locally Lipschitz vector field
\begin{equation}\label{eq:rkchecked-flow}
 J'_{ij}=-\sum_k w_k\frac{J_{ij}^2}{J_{ik}J_{jk}}.
\end{equation}

We show that this flow preserves $\mathsf P$. Write $\mathbf1$ for
the vector with every entry equal to one, and use $\succeq0$ and
$\preceq0$ for positive and negative semidefiniteness. At a matrix
$J$ with strictly positive entries and $\mathsf P(J)$, fix
$v\in\ker J$. For each $k$, let
\[
 D_k=\operatorname{diag}(J_{ik}),\quad
 H_k=D_k^{-1}JD_k^{-1},\quad
 U_k=J_{kk}^{-1}\mathbf1\mathbf1^T-H_k.
\]
Every entry in the $k$th row of $H_k$ equals $J_{kk}^{-1}>0$.
The Schur complement
at this diagonal entry is nonpositive, hence $U_k\succeq0$.
For $y_k=D_kv$ we have
$\mathbf1^Ty_k=0$, $H_ky_k=0$, and $U_ky_k=0$.
Expanding $H_k\odot H_k$, where $\odot$ denotes entrywise product,
therefore gives
\begin{equation}\label{eq:rkchecked-flow-tangent}
 v^TJ'v=-\sum_k w_k\,y_k^T(U_k\odot U_k)y_k\le0.
\end{equation}
The last inequality is the Schur product theorem.

We next justify forward invariance when several eigenvalues vanish.
For a one-by-one matrix, $\mathsf P$ always holds. Otherwise let
$a>0$ be the last diagonal entry, $b$ the remaining part of its
column, and $Q$ the Schur complement, so that
\[
 J=\begin{pmatrix}Q+bb^T/a&b\\b^T&a\end{pmatrix}.
\]
Then $\mathsf P(J)$ is equivalent to $Q\preceq0$. In the smooth
coordinates $(Q,b,a)$, write $Q'=\mathsf F(Q,b,a)$ for the locally Lipschitz
field induced by \eqref{eq:rkchecked-flow}. If $Q\preceq0$ and
$x\in\ker Q$, differentiation of the Schur complement and
\eqref{eq:rkchecked-flow-tangent}, applied to
$v=(x,-b^Tx/a)$, give $x^T\mathsf F(Q,b,a)x\le0$.

Let $Q_+$ be the positive semidefinite part of $Q$, obtained by
replacing each negative eigenvalue by zero, and put $Q_-=Q-Q_+$.
Near any point with $Q\preceq0$, the projected coordinates
$(Q_-,b,a)$ still represent a matrix with strictly positive entries.
For a real matrix $R$, write $\|R\|_F^2=\operatorname{tr}(R^TR)$
for its squared Frobenius norm. With a local Lipschitz constant $L_*$,
the squared distance
to the negative semidefinite cone satisfies
\[
 \frac12\frac{d}{dt}\|Q_+\|_F^2
 =\operatorname{tr}(Q_+\mathsf F(Q,b,a))
 \le L_*\|Q_+\|_F^2.
\]
Indeed, the term with $\mathsf F(Q_-,b,a)$ is nonpositive, because
$\operatorname{ran}Q_+\subset\ker Q_-$; the remaining term is
bounded by local Lipschitz continuity. The derivative formula
holds also at zero eigenvalues since $(\max\{x,0\})^2$ is $C^1$.
Gronwall's inequality proves local forward invariance. The entries of $J(t)$
remain strictly positive on the compact time interval, so a
first-exit argument gives invariance on $[0,1]$. Apply it at $t=1$.
\end{proof}

\subsection{Quadrature, regularization, and the boundary limit}

We prove Theorem~\ref{thm:rkchecked-bounded}. Green sections lie in
$L^2(q)$, and Poisson sections are uniformly bounded on the compact
integration support. The nonnegative resolvent series consequently gives
\begin{align*}
 M_q(x,y)&=\mathcal G_0(x,y)+
 \langle\mathcal G_0(x,\cdot),S\mathcal G_0(y,\cdot)\rangle_{L^2(q)},\quad x\ne y,\\
 K_q(z_0,s)&=P_s(z_0)+
 \langle\mathcal G_0(z_0,\cdot),SP_s\rangle_{L^2(q)},\\
 N_q(s,t)&=N_0(s,t)+\langle P_s,SP_t\rangle_{L^2(q)},\quad s\ne t.
\end{align*}
In particular $K_q(z_0,\cdot)$ is bounded and the off-diagonal
values of $N_q$ are finite. Set $M_q(x,x)=+\infty$.

Fix $0<\varepsilon\le1$. The operator $T_\varepsilon$ with kernel
$C_\varepsilon$ on $L^2(q)$ satisfies
$|T_\varepsilon f|\le T|f|$, hence $\|T_\varepsilon\|\le d$.
Partition a compact set containing the support into cells $E_i$ of
maximal diameter tending to zero. In each cell with
$w_i=\int_{E_i}q\,dA>0$, choose a node $z_i$; omit cells of zero mass.
Uniform continuity
of $C_\varepsilon$ makes the step operators converge in norm to
$T_\varepsilon$: their norm error is at most $\int q$ times the
uniform kernel error. Their nonzero blocks on normalized cell
indicators are
$\mathsf B_{ij}=\sqrt{w_i}C_\varepsilon(z_i,z_j)\sqrt{w_j}$.
For sufficiently fine partitions $\|\mathsf B\|\le(1+d)/2<1$.

Adjoin any fixed finite family of evaluation points, merging
coincident points and giving evaluation-only points weight zero.
Let $C$ be the matrix of $C_\varepsilon$ on this union and $W$ the
diagonal weight matrix. Then $W^{1/2}CW^{1/2}$ is $\mathsf B$ with
additional zero rows and columns. Applying
Lemma~\ref{lem:rkchecked-update} to $C,W$ proves $\mathsf P$ for
the reciprocal of the resulting resolvent and therefore for its
evaluation submatrix. Endpoint sections converge in $L^2(q)$.
The inverse identity and the uniform resolvent bound $2/(1-d)$
then show that the entries converge to
\[
 M_{q,\varepsilon}(x,y)=C_\varepsilon(x,y)+
 \langle C_\varepsilon(x,\cdot),
       (I-T_\varepsilon)^{-1}C_\varepsilon(y,\cdot)\rangle_{L^2(q)}.
\]
Thus its reciprocal evaluation matrices have $\mathsf P$.

Next let $\varepsilon\downarrow0$. Monotone convergence in the
nonnegative chain integrals and then in their series gives
$M_{q,\varepsilon}(x,y)\uparrow M_q(x,y)$ off the diagonal.
The limiting correction is bounded by
$(1-d)^{-1}\|\mathcal G_0(x,\cdot)\|_{L^2(q)}
\|\mathcal G_0(y,\cdot)\|_{L^2(q)}$.
On the diagonal the direct term tends to infinity. The finite
reciprocal matrices of $M_q$, with zero diagonal, therefore have
$\mathsf P$. When $q=0$, this conclusion follows directly from
$C_\varepsilon$ and requires no quadrature.

For distinct boundary points $s_1,\ldots,s_m$ and
$r>|z_0|$ tending to one, the disk Green formula gives
\[
 \frac{\mathcal G_0(rs,x)}{1-r}\longrightarrow P_s(x),\qquad
 \frac{\mathcal G_0(rs,rt)}{(1-r)^2}\longrightarrow N_0(s,t)
 \quad(s\ne t).
\]
The first convergence is uniform on the compact integration support.
The resolvent formulas hence give the corresponding limits
$M_q(z_0,rs)/(1-r)\to K_q(z_0,s)$ and
$M_q(rs,rt)/(1-r)^2\to N_q(s,t)$.
Congruence of the reciprocal matrix at
$rs_1,\ldots,rs_m,z_0$ by
$\operatorname{diag}(1-r,\ldots,1-r,1)$ and passage to the limit yield
\begin{equation}\label{eq:rkchecked-border}
 \mathsf P\!\left(\begin{pmatrix}B&\ell\\\ell^T&0\end{pmatrix}\right),
 \quad B_{ii}=0,\quad B_{ij}=N_q(s_i,s_j)^{-1}\ (i\ne j),\quad
 \ell_i=K_q(z_0,s_i)^{-1}.
\end{equation}

For any symmetric bordered matrix in \eqref{eq:rkchecked-border},
$x^TBx\le0$ on $\ell^\perp$. To see this, if $x^TBx>0$ there,
choose $y$ with $\ell^Ty\ne0$ and subtract a multiple of $x$ so
that $x^TBy=0$. A suitable last coordinate $c$ makes the quadratic
form positive on $(y,c)$ as well as $(x,0)$; these vectors are
orthogonal for that form and span a positive two-plane, a
contradiction. Taking $x_i=K_q(z_0,s_i)a_i$ proves conditional
negativity. Repeated points are combined first. This proves the theorem.

\subsection{A single truncation and Hilbert realization}

\begin{corollary}\label{cor:rkchecked-extension}
Let $q\ge0$ be measurable and satisfy
\eqref{eq:rkchecked-subcritical}. Fix $z_0\in\mathbb D$ and define
$K_q,N_q$ by their nonnegative chains. If $K_q(z_0,s)<\infty$ for
every $s$, and $N_q(s,t)<\infty$ for $s\ne t$, then $D_q$ in
\eqref{eq:rkchecked-D} is conditionally negative semidefinite.
In particular this holds for $q\in L^p(\mathbb D)$, $p>1$, by
Theorem~\ref{thm:hk-estimates}, whose proof in Sections 5.1--5.3
does not use reciprocal positivity.
\end{corollary}
\begin{proof}
For integers $m\ge2$, put
\[
 q_m(z)=\min\{q(z),m\}\mathbf1_{\{|z|\le1-1/m\}}.
\]
The chain formulas for the original $q$ are interpreted as
nonnegative kernel integrals; no global $L^2(q)$ representation is used.
These bounded compactly supported potentials increase to $q$ and
satisfy the same form bound. The theorem applies to each $q_m$.
Tonelli and monotone convergence, first for each chain and then for
the series, give $K_{q_m}\uparrow K_q$ and $N_{q_m}\uparrow N_q$
off the diagonal. Their finite positive limits imply
$D_{q_m}(s,t)\to D_q(s,t)$; monotonicity of this quotient is not
needed. Pass to the limit in each finite conditional-negative
inequality. Thus the limits have been taken in the order quadrature,
regularization, boundary approach, and finally $m\to\infty$.
\end{proof}

\begin{proposition}[Hilbert space realization]\label{prop:rkchecked-realization}
Under the hypotheses of Corollary~\ref{cor:rkchecked-extension},
assume also that $K_q(z_0,\cdot)$ is bounded. For any prescribed
$s_*\in\mathbb S^1$, there is an injective Lipschitz map
$b:\mathbb S^1\to\mathcal H$ into a real separable Hilbert space
such that $D_q(s,t)=\|b(s)-b(t)\|^2$ and $b(s_*)=0$.
\end{proposition}
\begin{proof}
Schoenberg's Gram construction \cite{Schoenberg1938} uses
\[
 H(s,t)=\tfrac12\bigl(D_q(s,s_*)+D_q(t,s_*)-D_q(s,t)\bigr).
\]
Appending the negative sum of any coefficient family at $s_*$
shows that $H$ is positive semidefinite. Complete the formal span
of the boundary points modulo the nullspace of this Gram form.
The resulting vectors $b(s)$ have the stated squared distances.
Since $N_q\ge N_0$,
\[
 \|b(s)-b(t)\|^2\le\pi\|K_q(z_0,\cdot)\|_\infty^2|s-t|^2.
\]
This proves Lipschitz continuity. Positivity of $K_q$ and finite
positive off-diagonal values of $N_q$ give injectivity. Restricting
to the closed span of the compact image gives separability.
\end{proof}

At $q=0$ the normalization is
\[
 D_0(s,t)=\frac1{4\pi}
 \left|\frac{s-z_0}{1-\overline{z_0}s}
       -\frac{t-z_0}{1-\overline{z_0}t}\right|^2,
\]
so $D_0(s,t)=|s-t|^2/(4\pi)$ when $z_0=0$.
Conditional negativity is closed under pointwise convergence of
finite kernels, also for nonmonotone approximating potentials.
For $q\in L^p$, $p>1$, the $C^{1,\beta}$ regularity required in
the multiplier argument is supplied separately by
Proposition~\ref{prop:hk-curve}.

\section*{Acknowledgment of AI assistance}
OpenAI's ChatGPT was used in the mathematical development and preparation
of this manuscript, including the boundary-primitive estimates, the
Gaussian multiplier argument, the approximation of spectral clusters, and
the stream-function proof, boundary angle estimates, and geometric
consequences in
Sections~\ref{subsec:global-eigen-coordinates}--\ref{subsec:eigen-coordinate-convexity}.
It also assisted with proof checking, literature comparison, the
adaptation of the reciprocal-matrix flow in Appendix~\ref{sec:rkchecked},
and revision of the exposition and notation.
The author takes full responsibility for the mathematical correctness and
attribution of the manuscript.

\end{document}